\documentclass[reqno]{amsart}
\usepackage{CJK}
\usepackage{hyperref}
\usepackage{cite}
\usepackage{amsmath}
\usepackage{mathrsfs}
\usepackage{color}

\makeatletter
\@namedef{subjclassname@2020}{%
	\textup{2020} Mathematics Subject Classification}
\makeatother

\numberwithin{equation}{section}

\newtheorem{thm}{Theorem}[section]

\newtheorem{lem}[thm]{Lemma}

\newtheorem{example}{Example}[section]
\newtheorem{Definition}[thm]{Definition}

\theoremstyle{definition}

\newtheorem{remark}[thm]{Remark}
\allowdisplaybreaks

\newcommand{\cA}{{\mathcal{A}}}

\newcommand{\supp}{\text{supp}}

\DeclareMathOperator*{\essinf}{ess\,inf}

\title[Renormalized and entropy solutions]{Renormalized  and entropy solutions to the  general nonlinear parabolic equations in Musielak-Orlicz spaces}
\author{Ying Li,  Chao Zhang$^*$}
\address{Ying Li\newline
School of Mathematics, Harbin Institute of Technology, Harbin 150001, China\newline
\texttt{lymath@hit.edu.cn}}
\address{Chao Zhang\newline
 School of Mathematics and Institute for Advanced Study in Mathematics, Harbin Institute of Technology, Harbin 150001, China
\newline
\texttt{czhangmath@hit.edu.cn}}
\thanks{$^*$ Corresponding author.}

\begin{document}

\maketitle

\begin{abstract}
We study the well-posedness of solutions to the general nonlinear parabolic equations with merely integrable data in time-dependent Musielak-Orlicz spaces. With the help of a density argument, we establish the existence and uniqueness of both renormalized and entropy solutions. Moreover, we conclude that the entropy and renormalized solutions for this equation are equivalent. Our results cover a variety of problems, including those with Orlicz growth, variable exponents, and double-phase growth.
\end{abstract}

\section{Introduction}
Suppose that $\Omega$ is a bounded domain of $\mathbb{R}^{N}$ with Lipschitz boundary $\partial\Omega$, and $T$ is a positive number. Denote $\Omega_T=\Omega\times (0,T]$, $\Sigma=\partial\Omega\times(0,T)$. This manuscript aims to establish the existence and uniqueness of both renormalized and entropy solutions for a general parabolic equation with merely integrable data in time-dependent spaces. Namely, we consider the following parabolic problem
\begin{align}\label{eq:main}
   \left\{
  \begin{array}{ccc}
\partial_{t} u- \text{div}~ \cA(t,x,\nabla u)= f(t,x) \quad &\mbox{in}& \Omega_{T}:=(0,T] \times \Omega,\\
u(t,x) = 0\quad  &\mbox{on}& \Sigma:=(0,T]\times \partial \Omega ,\\
u(\cdot,0)=u_{0}(\cdot)  \quad &\mbox{in}& \Omega, \quad \quad \quad  \quad \quad \quad
  \end{array}
  \right.
\end{align}
where $f\in L^1(\Omega_{T})$ and $u_{0}\in L^{1}(\Omega)$. We assume the following hypotheses:
\begin{itemize}
    \item [(A1)]  $\cA(t,x,\xi):(0, T)\times \Omega\times \mathbb{R}^{N}\rightarrow \mathbb{R}^{N}$ is a Carath\'{e}odory  function (i.e., measurable in $t$ and $x$ for fixed $\xi$ and continuous with respect to $\xi$ for fixed  $(t,x)$).

\item [(A2)] $\cA(t,x,0)=0$ for almost every $(t,x)\in \Omega_T$ and there exist an $N$-function $M$ (see Definition \ref{Def-N}) and constants $c_1\in (0,1)$,  such that for all $\xi\in \mathbb{R}^N$ we have
    \[ c_1 \left( M(t,x,|\xi|)+M^{*}(t,x,|\cA(t,x,\xi)|)\right) \leq\cA(t,x,\xi)\cdot\xi,   \]
where $M^*$ is the complementary function of $M$ (see~Definition~\ref{defgonge}).

\item [(A3)] $ \left[ \cA(t,x,\xi)-\cA(t,x,\eta)\right]\cdot(\xi-\eta)> 0$ for all $\xi,\eta$ in $\mathbb{R}^{N}$, $\xi \neq \eta$ and a.e. ~~~~$(t,x)\in \Omega_T$.
\label{f3}
\end{itemize}

Motivated by fluids of nonstandard rheology, we focus on the general form of growth conditions for the leading term of the operator, which makes Musielak-Orlicz spaces a suitable function space for the considered problem. We do not assume any growth condition of doubling type on the function $M$. Instead, we impose a condition that balances the behavior of $M$ with respect to its variable, ensuring that smooth functions are modularly dense in the related Sobolev-type space.

Musielak-Orlicz spaces, which generalize Orlicz, variable exponent and double-phase spaces,  have received significant attention since the seminal contributions of  Marcellini \cite{marcellini1991regularity,marcellini1989regularity}. There is a large amount of literature on PDE problems in the framework of Musielak-Orlicz spaces~\cite{colombo2015bounded,diening2011lebesgue,harjulehto2019orlicz}. We refer to \cite{F12, LZ} for the existence of weak solutions in isotropic, separable and reflexive Musielak-Orlicz-Sobolev spaces. Donaldson was the first to study nonlinear elliptic boundary value problems in non-reflexive Musielak-Orlicz Sobolev spaces, see~\cite{donaldson1971nonlinear}. Subsequent advancements were made by  Gossez \cite{gossez1974nonlinear, gossez1979orlicz}. The groundwork for studying PDE problems in anisotropic Musielak-Orlicz spaces was established in the works of  \cite{cianchi2000fully, cianchi2007symmetrization, klimov1976imbedding}. For more recent results concerning PDEs in Musielak-Orlicz spaces, the readers may refer to monograph \cite{BOOK} and the review paper \cite{chlebicka2018pocket}.

As we consider problems with data of low integrability,  it is reasonable to work with renormalized solutions and entropy solutions, as they require less regularity in the data than standard weak solutions. The concept of renormalized solutions was initially introduced by DiPerna and Lions in \cite{diperna1989cauchy} for analyzing the Boltzmann equation. At the same time,  B\'{e}nilan et al. proposed the notion of entropy solutions in \cite{Benilan95an} for the study of nonlinear elliptic problems. The existence of renormalized solutions within the context of variable exponents has been explored in~\cite{Bendahmane10,wittbold2010existence}.  We refer to \cite{nassar2014renormalized,redwane2010existence} for this issue in the nonreflexive Orlicz-Sobolev space. Also, there have been a large number of papers devoted to the study of renormalized solutions
for the partial differential equations in Musielak-Orlicz space, see \cite{chlebicka2019parabolic,  Chlebicka18jde,gwiazda2015renormalized, Gwiazda18jde,gwiazda2012corrigendum,gwiazda2012renormalized}.   Entropy solutions for the $p(x)$-Laplace equation were investigated in ~\cite{Sanc09entropy}, with further research on entropy solutions exhibiting Orlicz growth available in~\cite{zhang10entropy,zhang12on}. For studies on the existence of entropy solutions in Musielak-Orlicz spaces, see~\cite{Benkirane22strongly, Elarabi18entropy}.
Moreover, Droniou and Prignet established the equivalence between entropy solutions and renormalized solutions for parabolic problems with polynomial growth in~\cite{Droniou07}, which was subsequently studied for variable exponent and Orlicz growth in \cite{Zhang10en, zhang2010, zhang2017}, and for $(p(x),q(x))$ growth in parabolic equations in ~\cite{ABCJMAA}. Additionally, Li et al. \cite{li2021entropy} demonstrated this equivalence for elliptic equations with general growth in Musielak-Orlicz spaces.

Recently, Chlebicka, Gwiazda and Zatorska-Goldstein \cite{Chlebicka19re} showed the existence and uniqueness of renormalized solutions of the following equation
\begin{align*}
   \left\{
  \begin{array}{ccc}
\partial_{t} u- \text{div}~ \mathcal{A}(t,x,\nabla u)= f(t,x) \quad &\mbox{in}& \Omega_{T},\\
u(t,x) = 0\quad  &\mbox{on}& \Sigma,\\
u(0,x)=u_{0}(x)  \quad &\mbox{in}& \Omega,
  \end{array}
  \right.
\end{align*}
where $f\in L^1(\Omega_T)$, $u_0\in L^1(\Omega)$, and $\mathcal{A}$ was assumed to be controlled by a $N$-function. The authors therein employed a delicate time approximation method to achieve smoothness in the time direction. Unlike the proof in \cite{Chlebicka19re}, we do not require separate approximations for time and space variables. We assume that the regularity of the modular function is strong enough to ensure the density of smooth functions in the related Sobolev-type space. In fact, this density can be guaranteed by the balance condition (B) in the isotropic Musielak-Orlicz space, see~Lemma~\ref{modular density}. Utilizing this density result, we provide direct proof of the existence and uniqueness of renormalized solutions to equation \eqref{eq:main}, demonstrate that the renormalized solution is also an entropy solution,  and establish the uniqueness of the entropy solution, thereby showing the equivalence between entropy and renormalized solutions. It is worth noting that when
$M(t,x,\xi)=|\xi|^{p(t,x)}$, the existence and uniqueness of entropy solutions for  equation \eqref{eq:main} was discussed by Bul\'{i}\v{c}ek  and  Wo\'{z}nicki in ~\cite{bu23ar}. The results of this work extend the results presented in \cite{bu23ar, Chlebicka19re,li2021entropy,zhang2010}.

 We  rely on the density of smooth functions in a relevant function space to study problem~\eqref{eq:main}, namely
\[\textbf{W}(\Omega_T):=\bigg\{u\in W^{1,x}_0L_M(\Omega_T)\cap L^2(\Omega_T),\partial_t u\in W^{-1,x}L_{M^{*}}(\Omega_T)+L^2(\Omega_T)\bigg\},\]
where $W^{1,x}_0 L_M(\Omega_T)$ and $W^{-1,x}L_{M^{*}}(\Omega_T)$ are defined in Section 2.

To ensure the density, one may assume the regularity of $M$. Note that the smooth functions are dense in $\textbf{W}(\Omega_T)$  in the modular topology if the following balance condition holds, see Lemma~\ref{modular density}.

{\bf Balance condition $\mathsf{(B)}$.} If there exists a function $\varrho:[0,\infty)\times \mathbb{R}^+ \to \mathbb{R}^+$ which is non-decreasing with respect to each of the variables such that for  $(t,x)\in \Omega_T$ and $(\tau,y)\in \Omega_{T}$, 
\[M(t,x,s)\leq \varrho(|t-\tau|+c|x-y|,s) M(\tau,y,s)\quad \text{with}\quad \limsup_{\epsilon\to 0^+}\varrho(\epsilon,\epsilon^{-N})<\infty.\]

We point out that the Balanced condition is only used to ensure the density of smooth functions in our proof.  In addition, throughout the paper, we assume that the $N$-functions $M(t,x,s)$ satisfies the following Y-condition.

{\bf $\mathsf{(Y)}$-condition.} A $N$-function $M$ is said to satisfy the $Y$-condition on a segment $[a,b]$ of the real line $\mathbb{R}$, if either
\begin{itemize}
	\item[\textbf{(Y1)}]: there exist $q_0\in \mathbb{R}^+$ and $1\leq i\leq N$ such that  $x_i\in [a,b]\to M(t,x,s)$ is increasing when $x_i\geq q_0$ and decreasing when $x_i\leq q_0$, or vise versa,
\end{itemize}
or 
\begin{itemize}
	\item[\textbf{(Y2)}]: there exists $1\leq i\leq N$ such that for all $ s\geq 0$, the partial function $x_i \in [a,b]\to M(t,x,s)$ is monotone on $[a,b]$.
\end{itemize}
Here $x_i$ stands for the $i^{th}$ component of $x\in \Omega$.

In Musielak-Orlicz spaces, the norm Poincar\'{e} inequality is no longer true in general.
We remark that the $\mathsf{(Y)}$-condition is only used as a sufficient condition to obtain the norm Poincar\'{e} inequality (see Lemma~\ref{lemAYNEW} below),  which is crucial in the proof of Lemma~\ref{lem:var}. This condition covers the assumption given by Maeda~\cite{mae} to provide the Poincar\'{e} integral form for variable exponents. See \cite{AYnew}  for more information about this condition.

In \cite{Chlebicka19re}, the authors used different approximation methods for spatial and time directions, making it unnecessary to assume  $\mathsf{(Y)}$-condition to obtain information about the $L_M$ norm of the solution itself. In our situation, we do not separate space and time approximations. We rely on the density argument of Lemma~\ref{modular density}, which shows that for any test function $\phi$, it is necessary that $\phi\in L_M(\Omega_{T})$ and $\nabla \phi\in L_M(\Omega_{T})$. However, when making a priori estimates for the solutions to the approximation problem, we can only obtain gradient information about the approximate solutions due to the growth condition (A2). Thus, we need to impose that $N$-function $M$ satisfies $\mathsf{(Y)}$-condition. Nevertheless, we provide a more straightforward approach that reveals the intrinsic connection between renormalized and entropy solutions for such equations. In our future work, we will explore more general conditions that ensure smooth functions are dense in the anisotropic Musielak-Orilcz space. We believe that our method, with slight modifications, will be applicable to anisotropic Musielak-Orlicz spaces.

Before we proceed to define the renormalized and entropy solution to \eqref{eq:main}. We first introduce the truncation operator $T_{k}(r)$ as follows:
\begin{equation*}
  T_{k}(r)=  \left\{\begin{array}{cl}r \quad &\text{if}~|r|\leq k,\\
k\frac{r}{|r|} \quad &\text{if}~|r|> k.
\end{array}\right.
\end{equation*}
Its primitive $\Theta_{k}:\mathbb{R}\rightarrow \mathbb{R}^{+}$ defined by
\begin{align}
 \Theta_{k}(r):= \int^{r}_{0}T_{k}(s)ds =\left\{
  \begin{array}{ccc}
  &\!\!\!\!\!\!\frac{r^{2}}{2}             & \mbox{if}\quad~~|r|\leq k,\\
  &\!\!\!\!\!\! k|r|-\frac{k^{2}}{2}           & \mbox{if}\quad~~|r| > k.
  \end{array}
  \right.
\end{align}
It is obvious that $\Theta_{k}(r)\geq 0$ and $\Theta_{k}(r)\le k|r|$.

The definitions of entropy and renormalized solutions for problem (\ref{eq:main}) are as follows.
\begin{Definition}
    \label{def}
    A  function  $u \in C([0,T]; L^{1}(\Omega)) $ is  an entropy solution to  problem (\ref{eq:main}) if $u$ satisfies the following two conditions:
    \begin{itemize}
        \item [(E1)] $u$ is a measurable function satisfying 
        \[T_{k}(u)\in W_0^{1,x}L_M(\Omega_T)\quad \text{for each}~ k>0~\text{and}\quad \cA(t,x,\nabla T_{k}(u))\in L_{M^*}(\Omega_{T});\]
     
        \item [(E2)] For every $k>0$ and every $\phi\in C^{1}(\overline{\Omega}_{T})$ with $\phi|_{\Sigma}=0$, the inequality
\begin{eqnarray}\label{f5}
&&\int_{\Omega}\Theta_{k}(u-\phi)(T)\,dx -\int_{\Omega}\Theta_{k}(u_{0}-\phi(0))\,dx +\int^{T}_{0}\langle\phi_{t},T_{k}(u-\phi)\rangle \,dt \nonumber \\
& & \quad+\int^{T}_{0}\int_{\Omega}\cA(t,x,\nabla u)\cdot\nabla T_{k}(u-\phi)\,dxdt= \int^{T}_{0}\int_{\Omega}fT_{k}(u-\phi)\,dxdt
\end{eqnarray}
holds.
    \end{itemize}
\end{Definition}

\begin{Definition}
 A function $u\in  C([0,T];L^1(\Omega))$ is a renormalized solution to  problem (\ref{eq:main}) if $u$ satisfies the following two conditions:
    \begin{itemize}
        \item [(R1)]    $u$ is a measurable function satisfying 
        \[T_{k}(u)\in W_0^{1,x}L_M(\Omega_T)\quad \text{for each}~ k>0, \qquad \cA(t,x,\nabla T_{k}(u))\in L_{M^*}(\Omega_{T}) \]
        and 
$$\int_{\{l<|u|<l+1\}}\cA(t,x,\nabla u)\cdot \nabla u\,dxdt \rightarrow 0 \quad \quad \mbox{as}\quad  l  \rightarrow  +\infty;$$
\item [(R2)]   For every $ \phi \in C^{1}(\overline{\Omega}_{T})$ with $\phi(\cdot,T)=0$ and $S\in W^{2,\infty}(\mathbb{R})$ with
$S'$ has compact support, we have
\begin{equation}\label{R2}
    \begin{split}
     &-\int_{\Omega}S(u_0)\phi(x,0)\,dx-\int^{T}_{0}\int_{\Omega}S(u)\partial_{t}\phi\,dxdt\\
   &\quad +\int^{T}_{0}\int_{\Omega}\cA(t,x,\nabla u)\cdot \nabla\left(S'(u)\phi\right)\,dxdt =\int^{T}_{0}\int_{\Omega}fS'(u)\phi\,dxdt.
  \end{split}\end{equation}
    \end{itemize}

\end{Definition}

Now, let us state the main results of this work.

\begin{thm}\label{thm:renormalized}
{\it  Assume that $ f\in L^{1}(\Omega_{T})$, $u_{0}\in L^{1}(\Omega)$, $\cA$ satisfies the conditions $(A1)$--$(A3)$, $N$-function $M$ is regular enough so that the set of smooth functions is dense in $\textbf{W}(\Omega_T)$ in the modular topology, and $M$ satisfies $\mathsf{(Y)}$-condition. Then there exists a unique renormalized solution to  problem (\ref{eq:main}).}
\end{thm}

\begin{thm}\label{thm:entropy}
{\it  Assume that $ f\in L^{1}(\Omega_{T})$, $u_{0}\in L^{1}(\Omega)$, $\cA$ satisfies the conditions $(A1)$--$(A3)$,  $N$-function $M$ is regular enough so that the set of smooth functions is dense in $\textbf{W}(\Omega_T)$ in the modular topology, and $M$ satisfies $\mathsf{(Y)}$-condition.  Then the renormalized solution for problem \eqref{eq:main} is also an entropy solution to problem (\ref{eq:main}), and the entropy solution is unique.}
\end{thm}

\begin{remark}\label{916}
 {\it The entropy solution obtained in Theorem \ref{thm:entropy}} is equivalent to the renormalized solution of (\ref{eq:main}).
\end{remark}

\begin{example}  We will provide examples of time and space dependent Musielak-type spaces that are admissible in our investigation.
\begin{enumerate}
\item  If $M(t,x,s)=|s|^p$ with $1<p<\infty$, then the space $W^{1,x}L_M(\Omega_{T})$ is the classical Bochner space $L^p(0,T;W^{1,p}(\Omega))$ defined as 
\[ L^{p}(0,T;W^{1,p}(\Omega)):=\left\{u:(0,T)\to W^{1,p}(\Omega):\int^{T}_0\|u(t)\|^p_{W^{1,p}(\Omega)}\,dt<+\infty\right\}.\]

\item  If $M(t,x,s)=M(s)$, then we obtain the Orlicz-Sobolev space.

 \item If $M(t,x,s)=s^{p(t,x)}$, $1<p^-\leq p(t,x)\leq p^+<+\infty$,
 where $p(t,x)$ is log-H\"{o}lder continuous and $x_i\to p(t,x)$ is monotone on a compact subset of the real line $\mathbb{R}$,  then we obtain the variable Sobolev space. Here $x_i$ stands for the $i^{th}$ component of $x\in \Omega$.

  \item If $M(t,x,s)=s^p+a(t,x)s^q$, with  $0\leq a\in C^{0,\alpha}$, $\alpha\in (0,1]$, $1<p\leq q<+\infty$, $\frac{q}{p}\leq 1+\frac{\alpha}{n}$,  and $x_i\to a(t,x)$ is monotone on a compact subset of the real line $\mathbb{R}$, then we obtain the double-phase spaces. Here $x_i$ stands for the $i^{th}$ component of $x\in \Omega$.

\end{enumerate}

\end{example}

We organize this paper in the following framework. In Section 2, we state some basic results that will be used later. We will prove the main results in Section 3. In the following statement, $C$ stands for a constant, which may vary even within the same inequality.

\section{Functional setting and main tools}
In this section, we give some functional settings that will be used below. For a more complete discussion on this subject, we refer the readers to~\cite{BOOK, ACGA}.
We begin with the definition of $N$-functions.

 \begin{Definition}\label{Def-N}
 A function $M(\cdot,s):\Omega_T\times \mathbb{R}^+ \to \mathbb{R}^+$ is called an $N$-function if $M(\cdot,s)$ is a measurable function for every $s\geq0$, $M(t,x,\cdot)$ is strictly increasing with respect to last variable, and $M(t,x,\cdot)$ is a convex function for almost every $(t,x)\in \Omega_T$ with $M(t,x,0)=0$,  $M(t,x,s)\to +\infty$ as $s\to +\infty$ and 
\[\lim\limits_{s\to 0}\frac{M(t,x,s)}{s}=0 \quad \text{and}\quad \lim\limits_{s\to +\infty}\frac{M(t,x,s)}{s}=+\infty.\] 
\end{Definition}
\begin{Definition}\label{defgonge}
    The complementary function $M^*$ to an $N$-function $M$ in the sense of Young defined by
\begin{equation}\label{eqgonge}
    M^*(t,x, \xi_1):=\displaystyle\sup_{ \xi_2\geq 0}  \left[\xi_1\xi_2 -M(t,x,\xi_2)\right]
\end{equation}
for any~~$\xi_1 \geq 0 $ and a.e. $(t,x)\in \Omega_T$.
\end{Definition}

For an $N$-function,  we define the general Musielak-Orlicz class $\mathcal{L}_{M}(\Omega_T)$ as the set of all measurable functions $u(t,x):\Omega_T\to \mathbb{R}$ such that
$$\int^{T}_{0}\int_{\Omega}M(t,x,|u(t,x)|)\,dxdt < \infty.$$
The Musielak-Orlicz space $L_{M}(\Omega_T)$(resp. $E_M(\Omega_{T})$) is defined as the set of all measurable functions $u:\Omega_{T}\to \mathbb{R}$ such that 
\begin{eqnarray*}
	\int^T_0\int_{\Omega}M\left(t,x,\frac{|u(t,x)|}{\lambda}\right)\,dxdt<+\infty 
\end{eqnarray*} 
for some $\lambda>0$ (resp. for all $\lambda>0$). Equipped with Luxemburg  norm
\begin{equation*}
    \|u\|_{L_{M}(\Omega_T)}=\inf \bigg\{\lambda>0:\int^{T}_{0}\int_{\Omega}M\bigg(t,x,\frac{|u(t,x)|}{\lambda}\bigg)\,dxdt\leq 1\bigg\}.
\end{equation*}
Then $L_M(\Omega_{T})$ is a Banach space and  $E_M(\Omega_T)$ is its closed subset.

An $N$-function $M$ is  called locally integrable on $\Omega_T$, if for any constant number $c>0$ and for any compact $\Omega_T'$ of $\Omega_T$, the following holds
\[\int_{\Omega_T'}M(t,x,c)\,dxdt<+\infty.\]

We remark here that if an $N$-function $M$ satisfies the balanced condition (B), then the function is naturally locally integrable, see~\cite{MR4120278, MR4124438}. It is shown in \cite[Lemma 2.1]{ahmida2024} that the continuous embedding $L_M(\Omega_T)\hookrightarrow L^1(\Omega_T)$ holds if either $M^*$  is locally integrable, or $M$ satisfies $ \essinf\limits_{(t,x)\in\Omega_T} M(t,x,1)\geq c>0$.

\begin{Definition}
Suppose that the complementary $N$-function $M^*$ of $M$ is locally integrable on $\Omega_T$, we define
\begin{eqnarray}
    W^{1,x}L_{M}(\Omega_T):=\{u:\Omega_T\to \mathbb{R}:u\in L_{M}(\Omega_T),|\nabla u|\in L_{M}(\Omega_T)\}
\end{eqnarray}
and
\begin{eqnarray}
    W^{1,x}E_{M}(\Omega_T):=\{u:\Omega_T\to \mathbb{R}:u\in E_{M}(\Omega_T),|\nabla u|\in E_{M}(\Omega_T)\}.
\end{eqnarray}
We denote $\nabla u$ the vector gradient with respect to the space variable. These spaces are normed  by $\|u\|_{W^{1,x}L_{M}(\Omega_T)}:=\|u\|_{L_M(\Omega_T)}+\|\nabla u\|_{L_M(\Omega_T)}$, and then $W^{1,x}L_{M}(\Omega_{T})$ is a Banach space.
\end{Definition}

Let $X$ and $Y$ be subsets of $L^1(\Omega_T)$ not necessarily related by duality. We say $f_n \to f$ for $\sigma(X,Y)$ if
\begin{equation*}
	\int^{T}_0\int_\Omega f_n g \, dxdt \xrightarrow{n\to +\infty}\int^{T}_0 \int_\Omega fg \, dxdt
\end{equation*}
for all $g \in Y$. If $X=L_M(\Omega_T)$ and $Y=E_{M^\ast}(\Omega_T)$, we recover the weak-$\ast$ convergence and can also denote $f_n \stackrel{\ast}\rightharpoonup f$.

We define $W^1L_M(\Omega)$ (resp. $W^1E_M(\Omega)$) as the set of all measurable function $u:\Omega\to \mathbb{R}$, such that for all $|\alpha|\leq 1$, the function $|D^{\alpha} u|$ belong to $L_{M}(\Omega)$ (resp. $E_{M}(\Omega)$), that is, 
\[\int_{\Omega}M\left( x,\frac{|D^\alpha u|}{\lambda}\right)\,dx< +\infty\quad \text{for some $\lambda>0$ (resp. for all $\lambda>0$)}.\]

Note that if a $N$-function $M$ is locally integrable, then the set of $C_0^{\infty}(\Omega)$-functions is contained in $W^{1}E_{M}(\Omega)$. Therefore, the norm closure of $C_0^{\infty}(\Omega)$-functions in $W^{1}E_{M}(\Omega)$, denoted by $W_{0}^{1}E_{M}(\Omega)$, is well defined. Moreover, if the pair of complementary $N$-functions $(M, M^*) $  are  both locally integrable, then the space $ W_{0}^{1}L_{M}(\Omega)$, defined as the closure of $C_0^{\infty}(\Omega)$-functions with respect to the weak-* topology \( \sigma(L_{M}, E_{M^*}) \), is also well defined.

\begin{Definition}
    Suppose that the complementary $N$-function $M^*$ of $M$ is locally integrable on $\Omega_T$, we define
\begin{eqnarray}
    W_0^{1,x}L_{M}(\Omega_T):=\{u:(0,T)\to W_0^{1}L_M(\Omega):u\in L_{M}(\Omega_T),|\nabla u|\in L_{M}(\Omega_T)\}
\end{eqnarray}
and
\begin{eqnarray}
    W_0^{1,x}E_{M}(\Omega_T):=\{u:(0,T)\to W_0^{1}E_M(\Omega):u\in E_{M}(\Omega_T),|\nabla u|\in E_{M}(\Omega_T)\}.
\end{eqnarray}
These spaces are equipped with the norm $\|u\|_{W^{1,x}L_{M}(\Omega_T)}$.
\end{Definition}

\begin{lem}[Theorem 1.1, \cite{AYnew}]\label{lemAYNEW}
	Assume that the pair of complementary $N$-functions  $M$ and $M^*$ satisfy both balance condition $(B)$, and $M$ satisfies $\mathsf{(Y)}$-condition. Then there exists a constant $C$ depending only on $\Omega_{T}$ such that for every $u\in W^{1,x}_{0}L_M(\Omega_{T})$ it holds
	\[\|u\|_{L_M(\Omega_{T})}\leq C\|\nabla u\|_{L_M(\Omega_{T})}.\]
\end{lem}
From Lemma~\ref{lemAYNEW}, the two norms $\|\cdot\|_{L_M(\Omega_T)}+\|\nabla \cdot\|_{L_M(\Omega_T)}$ and $\|\nabla \cdot\|_{L_M(\Omega_T)}$ are equivalent on $W^{1,x}_0 L_{M}(\Omega_{T})$. In addition, it follows from  Definition~\ref{defgonge} that
\begin{equation}
    \int^{T}_{0}\int_{\Omega}|u_1 u_2|\,dxdt\leq 2\|u_1\|_{L_{M}(\Omega_T)}\|u_2\|_{L_{M^{*}}(\Omega_T)}
\end{equation}
for all $u_1\in L_{M}(\Omega_T)$ and $u_2\in L_{M^*}(\Omega_T)$.
We say that a sequence $\{v_n\}_{n=1}^{\infty}$ converges modularly to $v$ in $L_M(\Omega_T)$, if there exists $\lambda>0$ such that
$$\int^{T}_{0}\int_{\Omega}M\bigg(t,x,\frac{|v_n-v|}{\lambda}\bigg)\,dxdt\to 0 \ \ \ \text{as} \ \ \ n\to +\infty.$$
For the notion of this convergence, we write $v_n \xrightarrow[]{M}v$.

 \begin{lem}[Theorem 5.2,~\cite{ahmida2024}]\label{modular density}
Assume that $\cA$ satisfies the conditions $(A1)$--$(A3)$, and the pair of complementary $N$-functions  $M$ and $M^*$ satisfy both balance condition $(B)$. Then, for every $\phi \in \textbf{W}(\Omega_T)$ there exists a sequence $\{\phi_{\delta}\}_{\delta}\subset C_0^{\infty}((0,T];C^{\infty}_{0}(\Omega))$ such that
\begin{equation*}
 \begin{split}
 \phi_{\delta}\xrightarrow[]{}\phi \ \ \ \ \ \  \ \text{in} \ \ &L^2(\Omega_T), \\
 \partial_t \phi_{\delta}\xrightarrow[]{M}\partial_t\phi \ \ \ \ \text{in} \ \ \ &W^{-1,x}L_{M^*}(\Omega_T)+L^2(\Omega_T),\\
 D^{\alpha} \phi_{\delta}\xrightarrow[]{M} D^{\alpha} \phi, |\alpha|\leq 1 \ \ \  \text{in} \ \ &L_M(\Omega_T),\\
 \end{split}
\end{equation*}
where $W^{-1,x}L_{M^*}(\Omega_T)$ is defined as
\begin{equation*}
\begin{split}
W^{-1,x}L_{M^*}(\Omega_T):&=\left\{u:(0,T)\to W^{-1}L_{M^*}(\Omega):u=\tilde{u}-\text{\rm div}U,\right. \\
&\quad\left.~\textrm{with}~\tilde{u}\in L_{M^*}(\Omega_T)~\text{and}~U \in L_{M^*}(\Omega_T)\right\}.
\end{split}
\end{equation*}
\end{lem}

  The following fact is a consequence  of modular topology.
\begin{lem}[Lemma 2,~\cite{ACGA}]\label{lem:mo}
	Let $M$ be an $N$-function and  $u_n, u\in L_M(\Omega_{T})$. If $u_n\xrightarrow[]{M} u$ modularly, then $u_n\to u$ in $\sigma(L_M,L_{M^*})$.
\end{lem}

Using a proof strategy analogous to that in \cite{Elmahi05para}, we can establish the following lemma.
\begin{lem}\label{lem:var}
Suppose that $f\in C_0^\infty(\Omega_T)$ and $u_0\in C_0^\infty(\Omega)$, $\cA$ satisfies the conditions $(A1)$--$(A3)$, $N$-function $M$ is regular enough so that the set of smooth functions is dense in $\textbf{W}(\Omega_T)$ in the modular topology, and $M$ satisfies $\mathsf{(Y)}$-condition. Then there exists at least one distributional solution $u\in\textbf{W}(\Omega_T)$ of problem \eqref{eq:main} satisfying $u(x,0)=u_0(x)$ for almost every $x\in \Omega$. Furthermore, for all $\tau\in (0,T]$, we have
\begin{equation*}
-\int^{\tau}_0\int_{\Omega}\partial_t \phi u\,dxdt+\int_{\Omega}u\phi\,dx \bigg|^{\tau}_0+\int^{\tau}_0\int_{\Omega}\cA(t,x,\nabla u)\cdot \nabla \phi \,dxdt =\int^{\tau}_0\int_{\Omega}f\phi\, dxdt
\end{equation*}
for every $\phi \in \textbf{W}(\Omega_T)$.
\end{lem}

Then, we will present some preliminary lemmas that will be used later.

\begin{lem}[Lemma 2.7, \cite{zhang10entropy}]\label{w}{\it Let $\Omega_T$ be a measurable with finite Lebesgue measure, and let ${f_{n}}$ be a sequence of functions in $L^{p}(\Omega_T) (p \ge 1)$ such that
\begin{eqnarray*}
&& f_{n} \rightharpoonup f ~~~~~\quad ~\quad ~~\mbox{weakly in }~~~~~~~L^{p}(\Omega_T),\\
&& f_{n} \rightarrow g ~~~~~~~~~~\quad~~\quad~~~~\mbox{a.e.~ in }~~~~~~~~~~\Omega_T.
\end{eqnarray*}
Then $f=g$ a.e. in $\Omega_T$.
}
\end{lem}

\begin{lem}[Theorem 13.47, \cite{RKR}]\label{zui}
	Let $f_n, f\in L^1(\Omega_T)$ such that $ f_{n}\geq 0$ a.e. in $\Omega_T$, $f_{n} \rightarrow f$ a.e. in $\Omega_{T}$ and
	\begin{eqnarray*}
		 \int^{T}_0\int_{\Omega} f_n\,dxdt \to \int^{T}_0\int_{\Omega}f\,dxdt~\text{as}~n\to+\infty.
	\end{eqnarray*}
	Then $f_n\to f$ strongly in $L^1(\Omega_T)$.
\end{lem}

\section{The proofs of  main results}
We are now in a position to give the proofs of our main results. Some of the reasoning is based on the ideas developed in \cite{Chlebicka18jde, Gwiazda18jde,zhang2010}.

We first consider the following approximate problems
\begin{eqnarray}\label{eq:appro}
  \left\{\begin{array}{ccc}
\partial_t u_{n}- \text{div}~\cA(t,x,\nabla u_{n})= f_{n}  \quad  &\mbox{in}& \quad \Omega_{T},\\
u_{n} = 0           \quad  &\mbox{on}& \quad \Sigma,\\
u_{n}(0,x) = u_{0n} \quad   &\mbox{in}& \quad  \Omega,
  \end{array}
  \right.
\end{eqnarray}
where the two sequences of functions $\{f_{n}\}\subset C^{\infty}_{0}(\Omega_{T})$ and $\{u_{0n}\}\subset C^{\infty}_{0}(\Omega)$ strongly convergent respectively to $f$ in $L^{1}(\Omega_{T})$ and to $u_{0}$ in $L^{1}(\Omega)$ such that
\begin{eqnarray}\label{3.1}
\|f_{n}\|_{L^{1}(\Omega_{T})} \leq \|f\|_{L^{1}(\Omega_{T})},\quad \quad ~~~~~ \|u_{0n}\|_{L^{1}(\Omega)}\leq \|u_{0}\|_{L^{1}(\Omega)}.
\end{eqnarray}

It follows from Lemma~\ref{lem:var} that there exists a distributional  solution $u_n\in\textbf{W}(\Omega_T)$ for problem~\eqref{eq:appro}, such that
 \begin{equation}\label{eq:appro1}
\int^{\tau}_0\int_{\Omega}\partial_tu_n\phi\,dxdt+\int^{\tau}_0\int_{\Omega}\cA(t,x,\nabla u_n)\cdot \nabla \phi \,dxdt =\int^{\tau}_0\int_{\Omega}f_n\phi\, dxdt
\end{equation}
for every $\phi \in \textbf{W}(\Omega_T)$.

Taking the test function  as $T_{k}(u_{n})\chi_{(0,\tau)}$ with $\tau \in (0,T]$ in \eqref{eq:appro1}, we have
\begin{equation}\label{eq:s1}
\begin{split}
&\int_{\Omega}\Theta_{k}(u_{n})(\tau)\,dx -\int_{\Omega}\Theta_{k}(u_{0n})\,dx \\
 & \quad+\int^{\tau}_{0}\int_{\Omega}\cA(t,x,\nabla T_{k}(u_{n}))\cdot \nabla T_{k}(u_{n})\,dxdt = \int_{0}^{\tau}\int_{\Omega}f_{n}T_{k}(u_{n})\,dxdt.
\end{split}\end{equation}
According to the definition of $\Theta_{k}(r)$ and (\ref{3.1}), we deduce
\begin{equation}\label{3.4}
\begin{split}
&\int^{\tau}_{0}\int_{\Omega}\cA(t,x,\nabla T_{k}(u_{n}))\cdot \nabla T_{k}(u_{n})\,dxdt + \int_{\Omega}\Theta_{k}(u_{n})(\tau)\,dx\\
&\quad \leq k \left( \|f_{n}\|_{L^{1}(\Omega_{T})}+\|u_{0n}\|_{L^{1}(\Omega)} \right) \leq k\left( \|f\|_{L^{1}(\Omega_{T})}+\|u_{0}\|_{L^{1}(\Omega)}\right).
\end{split}\end{equation}
Recalling condition (A2), we have
\begin{eqnarray}\label{3.5}
c_1\int^{\tau}_{0}\int_{\Omega}M(t,x,|\nabla T_{k}(u_{n})|)\,dxdt \leq  \int^{\tau}_{0}\int_{\Omega}\cA(t,x,\nabla T_{k}(u_{n}))\cdot \nabla T_{k}(u_{n})\,dxdt\leq Ck
\end{eqnarray}
and
\begin{eqnarray}\label{2224}
	c_1\int^{\tau}_{0}\int_{\Omega}M^*(t,x,|\cA(t,x,\nabla T_{k}(u_{n})|)\,dxdt \leq Ck.
\end{eqnarray}
Since  $L_M(\Omega_T)\hookrightarrow L^1(\Omega_T)$, we know 
\begin{eqnarray}\label{3.6}
\int^{\tau}_{0}\int_{\Omega}|\nabla T_{k}(u_{n})|\,dxdt\le C\left(k+1\right),
\end{eqnarray}
that is,  $T_{k}(u_{n})$ is bounded in $L^{1}(0,T; W^{1,1}_{0}(\Omega))$.

 Choosing $k=1$ in the inequality (\ref{3.4}), we find  that
\begin{eqnarray*}
\int_{\Omega}\Theta_{1}(u_{n}(\tau))dx \leq \|f\|_{L^{1}(\Omega_{T})}+\|u_{0}\|_{L^{1}(\Omega)}
\end{eqnarray*}
for a.e. $\tau \in (0.T]$. Moreover,
\begin{eqnarray*}
\int_{\Omega}|u_{n}(\tau)|dx \le \text{meas}(\Omega)+\|f\|_{L^{1}(\Omega_{T})}+\|u_{0}\|_{L^{1}(\Omega)}.
\end{eqnarray*}
Therefore, we obtain
\begin{eqnarray}\label{3.9}
\|u_{n}\|_{L^{\infty}(0,T;L^{1}(\Omega))}\le C.
\end{eqnarray}

\vskip .2cm

 Now, we are ready to prove the existence and uniqueness of renormalized solutions.

\vskip .2cm

\noindent\textbf{Proof of Theorem \ref{thm:renormalized}}.
\textbf{(1). Existence of renormalized solutions.}

We divide our proof into several steps.

\noindent\textbf{Step 1}. Prove the convergence of  $\{u_n\}$ in $C([0,T]; L^1(\Omega))$ ans finds its subsequence which is almost everywhere convergent in $\Omega_T$.

Due to \eqref{eq:appro1}, we can write the weak form as
\begin{equation}\label{eq:almost}
    \begin{split}
    &\int^{T}_{0} \langle \partial_t(u_{n}-u_{m}),\phi\rangle\,dt +\int^{T}_{0}\int_{\Omega}[\cA(t,x,\nabla u_{n})-\cA(t,x,\nabla u_{m})]\cdot \nabla \phi\,dxdt \\
&\quad=\int^{T}_{0}\int_{\Omega}(f_{n}-f_{m})\phi \,dxdt
 \end{split}\end{equation}
for all $m,n \in \mathbb{Z}$ and $\phi\in  \textbf{W}(\Omega_T)$.
It follows from the Fenchel's Young inequality and \eqref{3.4} that
\begin{eqnarray*}
|\cA(t,x,\nabla u_{n})\cdot \nabla u_{m}|\le \left(M(t,x,|\nabla u_{n}|)+M^{*}\left(t,x,|\cA(t,x,\nabla u_{m})|\right)\right) \in L^{1}(\Omega_{T}).
\end{eqnarray*}
Define
\begin{eqnarray}\label{312}
\alpha_{n,m}:=\int^{T}_{0}\int_{\Omega}|f_{n}-f_{m}|\,dxdt +\int_{\Omega}|u_{0n}-u_{0m}|dx.
\end{eqnarray}
Since $\{f_{n}\}$ and $\{u_{0n}\}$ are convergent in $L^{1}$, we have
\begin{eqnarray*}
\lim_{n,m \to+ \infty}\alpha_{n,m}=0.
\end{eqnarray*}
Taking $w=T_{1}(u_{n}-u_{m})\chi_{(0,\tau)}$ with $\tau\le T$ as a test function in \eqref{eq:almost}, and discarding the positive term we obtain
\begin{eqnarray*}
\int_{\Omega}\Theta_{1}(u_{n}-u_{m})(\tau)dx&\le&\int_{\Omega}\Theta_{1}(u_{0n}-u_{0m})dx+\|f_{n}-f_{m}\|_{L^{1}(\Omega_{T})}\nonumber\\
&\le&\|u_{0n}-u_{0m}\|_{L^{1}(\Omega)}+\|f_{n}-f_{m}\|_{L^{1}(\Omega_{T})}=\alpha_{n,m}.
\end{eqnarray*}
Therefore, we conclude that
\begin{eqnarray*}
\int_{\{|u_{n}-u_{m}|<1\}}\frac{|u_{n}-u_{m}|^{2}(\tau)}{2}\,dx +\int_{\{|u_{n}-u_{m}|\ge1 \}}\frac{|u_{n}-u_{m}|(\tau)}{2}\,dx\nonumber\\
\le \int_{\Omega}[\Theta_{1}(u_{n}-u_{m})](\tau)\,dx \le \alpha_{n,m}.
\end{eqnarray*}
Moreover,
\begin{eqnarray*}
\int_{\Omega}|u_{n}-u_{m}|(\tau)\,dx &=&\int_{\{|u_{n}-u_{m}|<1\}}|u_{n}-u_{m}|(\tau)\,dx +\int_{\{|u_{n}-u_{m}|\ge 1\}}|u_{n}-u_{m}|(\tau)\,dx\nonumber\\
&\le &\left(\int_{\{|u_{n}-u_{m}|<1\}}|u_{n}-u_{m}|^{2}(\tau)\,dx\right)^{\frac{1}{2}}\text{meas}(\Omega)^{\frac{1}{2}}+2\alpha_{n,m}\nonumber\\
&\le &\left( 2\text{meas}(\Omega)\right)^{\frac{1}{2}}\alpha^{\frac{1}{2}}_{n,m}+2\alpha_{n,m}.
\end{eqnarray*}
Thus,  we deduce that
\begin{eqnarray*}
\|u_{n}-u_{m}\|_{C([0,T];L^{1}(\Omega))}\rightarrow 0 ~~~~~\quad \mbox{as}~ ~n,~m \rightarrow +\infty,
\end{eqnarray*}
 which implies that $\{u_{n}\}$ is a Cauchy sequence in $C([0,T];L^{1}(\Omega))$. Then $u_{n}$ converges to $u$ in $C([0,T];L^{1}(\Omega))$. We find an a.e. convergent subsequence (still denoted by $\{u_{n}\})$ in $\Omega_{T}$ such that
\begin{eqnarray}\label{3}
u_{n}\rightarrow u ~~~\quad \mbox{a.e.~~~in ~~~~~} ~~~~~\Omega_{T}\qquad ~~~~~~\mbox{as}~~~~~n\rightarrow +\infty.
\end{eqnarray}

\noindent\textbf{Step 2.} Show that the sequence $\{\nabla u_{n}\}$ converges almost everywhere in $\Omega_{T}$ to $\nabla u$ (up to a subsequence). Let us first set $\delta>0$ and denote
\begin{eqnarray*}
&&A_{1}:=\{ (t,x)\in \Omega_{T}:|\nabla u_{n}|>h\}\cup \{(t,x)\in \Omega_{T}:|\nabla u_{m}|>h\},\\
&&A_{2}:=\{ (t,x)\in \Omega_{T}:|u_{n}-u_{m}|>1\}
\end{eqnarray*}
and
\begin{eqnarray*}
A_{3}:=\{ (t,x)\in\Omega_{T}:|\nabla u_{n}|\le h,|\nabla u_{m}|\le h,|u_{n}-u_{m}|\le1,|\nabla u_{n}-\nabla u_{m}|>\delta\},
\end{eqnarray*}
where $h$ will be chosen later. Next, we shall show that  $\{\nabla u_{n}\}$ is a Cauchy sequence in measure.
It is easy to check that
\begin{eqnarray*}
\{ (t,x)\in\Omega_{T}:|\nabla u_{n}-\nabla u_{m}|>\delta \}\subset A_{1}\cup A_{2}\cup A_{3}.
\end{eqnarray*}
\noindent Firstly, we notice that
\begin{eqnarray*}
\{ (t,x)\in\Omega_{T}:|\nabla u_{n}|\ge h\}\subset \{ (t,x)\in\Omega_{T}:|u_{n}|\ge k\}\cup\{ (t,x)\in\Omega_{T}:|\nabla T_{k}(u_{n})|\ge h\}
\end{eqnarray*}
for all $k>0$.
Thus, using (\ref{3.9}) and (\ref{3.6}), we know that there exist constants $C>0$ such that
\begin{eqnarray*}
\mbox{meas}\{ (t,x)\in \Omega_{T}:|\nabla u_{n}|\ge h\} \le \frac{C}{k}+\frac{C(k+1)}{h},
\end{eqnarray*}
when $h$ is large appropriately. By choosing $k=Ch^{\frac{1}{2}}$, we deduce that
\begin{eqnarray*}
\mbox{meas}\{ (t,x)\in \Omega_{T}:|\nabla u_{n}|\ge h\} \le Ch^{-\frac{1}{2}}.
\end{eqnarray*}
Let $\varepsilon>0$. We may let $h=h(\varepsilon)$ large enough such that
\begin{eqnarray}\label{A1}
\mbox{meas}(A_{1})\le \frac{\varepsilon}{3} \quad \mbox{for~~ all} ~n,m>0.
\end{eqnarray}
\noindent Secondly, by Step 1 we know that $\{u_{n}\}$ is a Cauchy sequence in measure. Then there exists $N_{1}(\varepsilon)\in \mathbb{N}$ such that
\begin{eqnarray}\label{A2}
\mbox{meas}(A_{2}) \le \frac{\varepsilon}{3} \qquad \mbox{for all}~n,m\ge N_{1}(\varepsilon).
\end{eqnarray}
\noindent Finally, from condition $(A3)$, we know that there exists a real-valued function $m(h,\delta)>0$ such that
\begin{eqnarray*}
[\cA(t,x,\nabla \eta)-\cA(t,x,\nabla \zeta)]\cdot({\nabla\eta-\nabla \zeta})\ge m(h,\delta)>0
\end{eqnarray*}
for all $\eta,\zeta \in \mathbb{R}^{N}$ with $|\eta|,|\zeta|\le h,\delta \le|\eta-\zeta|$.
By taking $T_{1}(u_{n}-u_{m})$ as a test function in \eqref{eq:almost} and integrate on $A_{3}$, we obtain
\begin{equation*}
    \begin{split}
       m(h,\delta) \mbox{meas} (A_{3})&\le \int_{A_{3}}[\cA(t,x,\nabla u_{n})-\cA(t,x,\nabla u_{m})]\cdot(\nabla u_{n}-\nabla u_{m})\,dxdt\\
&\le\int^{T}_{0}\int_{\Omega}[\cA(t,x,\nabla u_{n})-\cA(t,x,\nabla u_{m})]\cdot(\nabla u_{n}-\nabla u_{m})\,dxdt\\
&\le \int^{T}_{0}\int_{\Omega}|f_{n}-f_{m}|\,dxdt+\int_{\Omega}|u_{0n}-u_{0m}|dx =\alpha(n,m),
    \end{split}
\end{equation*}
which implies that
\begin{eqnarray}\label{66}
\mbox{meas}(A_{3}) \le \frac{\alpha(n,m)}{m(h,\delta)}\le \frac{\varepsilon}{3} \qquad \mbox{for all} ~n,m \ge N_{2}(\varepsilon,\delta).
\end{eqnarray}
Combining the estimates (\ref{A1})--(\ref{66}), we obtain
\begin{eqnarray*}
\mbox{meas}\{ (t,x)\in \Omega_{T}:|\nabla u_{n}-\nabla u_{m}|> \delta \}\le \varepsilon  \qquad \mbox{for all}~~ n,m \ge \rm {max}\{N_{1},N_{2}\},
\end{eqnarray*}
that is $\{\nabla u_{n}\}$ is a Cauchy sequence in measure. Therefore, we obtain a subsequence of $\{\nabla u_{n}\}$ which is  almost everywhere convergent in $\Omega_{T}$.  Moreover, a priori estimate (\ref{3.5}) and weak lower semi-continuity of a convex functional give that
\begin{eqnarray}\label{eq:Tk}
T_{k}(u_{n}) \rightharpoonup  T_{k}(u) ~~~~~~~~~\quad \mbox{ weakly ~in} ~~~L^{1}(0,T; W^{1,1}_{0}(\Omega)).
\end{eqnarray}
Therefore, we deduce from Lemma \ref{w} that
\begin{eqnarray}\label{twsl}
\nabla u_{n}\rightarrow \nabla u  \qquad \mbox{a.e. in}~\Omega_{T} ~\mbox{as}~~~~~n\rightarrow +\infty .
\end{eqnarray}

\noindent\textbf{Step 3.} Prove a decay condition for $u_n$. In this step, we aim to prove that
\begin{equation}
\lim\limits_{l\to +\infty}\lim\limits_{n\to +\infty}\int_{\{l\leq |u_n|\leq l+1\}}\cA(t,x,\nabla u_n)\cdot \nabla u_n\,dxdt=0.
\end{equation}
Define the function $T_{l,a}(s)=T_{a}(s-T_{l}(s))$ as
\begin{align*}
 T_{l,a}(s):=\left\{
  \begin{array}{ccc}s-l\rm{sign}(s)             & \mbox{if}&~~ l \leq |s| < l+a,\\
   a\rm{sign}(s)& \mbox{if}&~~l+a \leq |s|,\\
  0& \mbox{if}&~~ |s|< l.
  \end{array}
  \right.
\end{align*}
 Using $T_{l,a}(u_{n})=T_{a}\left(u_{n}-T_{l}(u_{n})\right)$ as a test function in \eqref{eq:appro1}, we find
 \begin{equation*}
 \begin{split}
&\int_{\{|u_{n}|>l\}}\Theta_{a}(u_{n}\mp l)(T)\,dx -\int_{\{|u_{0n}|>l\}}\Theta_{a}(u_{0n}\mp l)\,dx\\
&\quad +\int_{\{l \leq |u_{n}|\leq l+a\}}\cA(t,x,\nabla u_{n})\cdot \nabla u_{n}\,dxdt\\
 &\leq \int^{T}_0\int_{\Omega}f_{n}T_{l,a}(u_{n})\,dxdt,
\end{split} \end{equation*}
which yields that
\begin{eqnarray}\label{eq:step31}
\int_{\{l\leq |u_{n}|\leq l+a\}} \cA(t,x,\nabla u_{n})\cdot \nabla u_{n}\,dxdt \le a (\int_{\{|u_{n}|>l\}}|f_{n}|\,dxdt+\int_{\{|u_{0n}|>l\}}|u_{0n}|\,dx).
 \end{eqnarray}
 Recalling the convergence of $\{u_{n}\}$ in $C\left([0,T];L^{1}(\Omega)\right)$, we have
\begin{eqnarray*}
\lim\limits_{l\rightarrow +\infty} \text{meas}\{ (t,x)\in \Omega_T:|u_{n}|>l\}=0 \quad \quad \mbox{uniformly  with respect to} ~n.
 \end{eqnarray*}
 Therefore, passing to the limit first in $n$ then in $l$, we conclude from \eqref{eq:step31} that
 \begin{eqnarray*}
 \lim_{l  \to+ \infty}\lim\limits_{n\rightarrow +\infty} \int_{\{ l \leq |u_n| \leq l+a\}} \cA(t,x,\nabla u_n)\cdot \nabla u_n \,dxdt=0.
 \end{eqnarray*}
 Choosing $a=1$, we have
  \begin{eqnarray}\label{eqstep3}
\lim_{l\to+\infty}\lim\limits_{n\rightarrow +\infty} \int_{\{ l \leq |u_n| \leq l+1\}} \cA(t,x,\nabla u_n)\cdot \nabla u_n \,dxdt=0.
 \end{eqnarray}

\noindent\textbf{Step 4.} Establish the convergence of $\cA(t,x,\nabla T_{k}(u_n))\cdot\nabla T_{k}(u_n)$.

We are going to show that
    \begin{eqnarray}\label{goal}
        \cA(t,x,\nabla T_{k}(u_n)) \nabla T_k(u_n)  \rightarrow \cA(t,x,\nabla T_{k}(u))\nabla T_{k}(u) ~~~~~\quad ~\quad ~~\mbox{strongly in }~~~~~~~L^{1}(\Omega_T).
    \end{eqnarray}
In fact, according to \eqref{2224} and condition (A2), we obtain
\begin{eqnarray*}
    \int^{T}_{0}\int_{\Omega}M^*(t,x,|\cA(t,x,\nabla T_{k}(u_n))|)\,dxdt \leq C.
\end{eqnarray*}
Therefore, there exists a subsequence  of $\cA(t,x,\nabla T_{k}(u_n)$  such that as $n\to +\infty$,
\begin{eqnarray}\label{eq:step41}
    \cA(t,x,\nabla T_{k}(u_n))  \rightharpoonup \cA_{k}\quad \quad \text{weakly-$*$~in}~ L_{M^*}(\Omega_T).
\end{eqnarray}
Recalling the fact that $\nabla u_n \to \nabla u$ a.e. in $\Omega_T$ and $\cA(t,x,\nabla T_{k}(u_n))$ is continuous with respect to $\nabla T_{k}(u_n)$, we deduced that
\begin{equation*}
    \cA_{k}=\cA(t,x,\nabla T_{k}(u)).
\end{equation*}

To prove \eqref{goal}, we first show
\begin{eqnarray*}
        \lim_{n\to+\infty}\int^{T}_{0}\int_{\Omega}\cA(t,x,\nabla T_{k}(u_n))\cdot\nabla T_{k}(u_n)\,dxdt=\int^{T}_{0}\int_{\Omega}\cA(t,x,\nabla
 T_{k}(u))\cdot\nabla T_{k}(u)\,dxdt.
\end{eqnarray*}

Since $T_{k}(u)\in \textbf{W}(\Omega_T)$, it follows from  Lemma~\ref{modular density} that there exists a sequence $\{T_{k}(u)\}_{\delta}\subset C_0^{\infty}((0,T];C^{\infty}_{0}(\Omega))$ such that
\begin{equation*}
 \begin{split}
 \{T_{k}(u)\}_{\delta}\xrightarrow[]{M}T_{k}(u)& \ \ \ \text{in} \ \ L_M(
 \Omega_T), \\
 \nabla \{T_{k}(u_n)\}_{\delta}\xrightarrow[]{M}\nabla T_{k}(u)& \ \ \ \text{in} \ \ L_M(
 \Omega_T),\\
 \partial_t \{T_{k}(u)\}_{\delta}\xrightarrow[]{M} \partial_t T_{k}(u)& \ \ \ \text{in}\ \ W^{-1,x}L_{M^*}(\Omega_T)+L^2(\Omega_T).
 \end{split}
\end{equation*}
Define $\psi_l(r)=\min\{(l+1-|r|)^+,1\}$. Taking $\psi_{l}(u_n)\{T_{k}(u)\}_{\delta}\chi_{(0,\tau)}$ as a test function for problem \eqref{eq:appro}, we have
\begin{equation}\label{eqnew}
    \begin{split}
     &\int^{\tau}_0\int_{\Omega}\partial_tu_n \psi_{l}(u_n)(T_{k}(u))_{\delta}\,dxdt +\int^{\tau}_0\int_{\Omega}\cA(t,x,u_{n})\cdot \nabla\left(\psi_{l}(u_n)(T_{k}(u))_{\delta}\right)\,dxdt \\
&\quad=\int^{\tau}_0\int_{\Omega}f_{n}\psi_{l}(u_n)(T_{k}(u))_{\delta}\,dxdt.
    \end{split}
\end{equation}

 For the first term on the left-hand side of equation~\eqref{eqnew}, recalling the fact that $u_n\to u$ a.e in $\Omega_T$ as $n\to +\infty$, we deduced from Lebesgue dominated convergence theorem that
 \begin{equation*}
 \begin{split}
&\lim\limits_{\delta\to 0}\lim\limits_{n\to +\infty}\int^{\tau}_0\int_{\Omega}\partial_tu_n \psi_{l}(u_n)(T_{k}(u))_{\delta}\,dxdt \\
&=\lim\limits_{\delta\to 0}\lim\limits_{n\to +\infty}\int^{\tau}_0\int_{\Omega} \partial_t \int^{u_n}_0\psi_{l}(r)\,dr (T_{k}(u))_{\delta}\,dxdt\\
&=\lim\limits_{\delta\to 0}\lim\limits_{n\to +\infty}\bigg[\int_{\Omega}\int^{u_n}_0\psi_{l}(r)\,dr (T_{k}(u))_{\delta}\,dx\bigg|^{\tau}_0-\int^\tau_0\int_{\Omega}\int^{u_n}_0\psi_{l}(r)\,dr \partial_t (T_{k}(u))_{\delta}\,dxdt\bigg]\\
&=\int^{\tau}_0\int_{\Omega}\partial_tu\psi_{l}(u)T_{k}(u)\,dxdt.
 \end{split}
 \end{equation*}
 Taking $l\to +\infty$, it is obvious that
 \begin{equation}\label{eqnew1}
 \lim\limits_{l\to +\infty}\lim\limits_{\delta\to 0}\lim\limits_{n\to +\infty}\int^{\tau}_0\int_{\Omega}\partial_tu_n \psi_{l}(u_n)(T_{k}(u))_{\delta}\,dxdt =\int^{\tau}_0\int_{\Omega}\partial_tu T_{k}(u)\,dxdt.
\end{equation}

For the second term on the left-hand side of equation~\eqref{eqnew}, we know that
\begin{equation*}
\begin{split}
&\int^{\tau}_0\int_{\Omega}\cA(t,x, \nabla u_{n})\cdot \nabla\left(\psi_{l}(u_n)(T_{k}(u))_{\delta}\right)\,dxdt\\
&\quad= \int^{\tau}_0\int_{\Omega}\cA(t,x, \nabla u_{n})\cdot \nabla T_{l+1}(u_n) \psi_l'(u_n)(T_{k}(u))_{\delta}\,dxdt\\
&\quad\quad +\int^{\tau}_0\int_{\Omega}\cA(t,x, \nabla u_{n})\cdot \nabla (T_{k}(u))_{\delta}  \psi_l(u_n)\,dxdt\\
&\quad:=Z_1+Z_2.
\end{split}
\end{equation*}
For $Z_1$, it follows from~\eqref{eqstep3} that
\begin{equation*}
\begin{split}
&\lim\limits_{l\to +\infty}\lim\limits_{\delta\to 0}\lim\limits_{n\to +\infty}\int^{\tau}_0\int_{\Omega}\cA(t,x, \nabla u_{n})\cdot \nabla T_{l+1}(u_n) \psi_l'(u_n)(T_{k}(u))_{\delta}\,dxdt\\
&\quad \leq \lim\limits_{l\to +\infty}\lim\limits_{n\to +\infty}C\int_{\{l\leq |u_n|\leq l+1\}}\cA(t,x, \nabla u_{n})\cdot \nabla T_{l+1}(u_n)\,dxdt\to 0.
\end{split}
\end{equation*}
For $Z_2$, we deduced from \eqref{eq:step41} that 
\begin{equation}
\cA(t,x, \nabla T_{l+1}(u_{n}))\psi_l(u_n)\rightharpoonup \cA(t,x, T_{l+1}(\nabla u))\psi_l(u) ~\text{weakly in}~L^1(\Omega_T)~\text{as}~n\to +\infty,
\end{equation}
which gives that 
\begin{equation*}
\begin{split}
&\lim\limits_{l\to +\infty}\lim\limits_{\delta\to 0}\lim\limits_{n\to +\infty}\int^{\tau}_0\int_{\Omega}\cA(t,x, \nabla u_{n})\cdot \nabla (T_{k}(u))_{\delta}  \psi_l(u_n)\,dxdt\\
&\quad=\lim\limits_{l\to +\infty}\lim\limits_{\delta\to 0}\lim\limits_{n\to +\infty}\int^{\tau}_0\int_{\Omega}\cA(t,x, \nabla T_{l+1}(u_{n}))\cdot \nabla (T_{k}(u))_{\delta}  \psi_l(u_n)\,dxdt\\
&\quad=\int^{\tau}_0\int_{\Omega}\cA(t,x, \nabla u)\cdot \nabla T_{k}(u)\,dxdt.
\end{split}
\end{equation*}
The limit as $\delta\to 0$ results from Lemma \ref{lem:mo}.
Therefore, we have
\begin{equation}\label{eqnew2}
\begin{split}
&\lim\limits_{l\to +\infty}\lim\limits_{\delta\to 0}\lim\limits_{n\to +\infty}\int^{\tau}_0\int_{\Omega}\cA(t,x, \nabla u_{n})\cdot \nabla\left(\psi_{l}(u_n)(T_{k}(u))_{\delta}\right)\,dxdt\\
&\quad=\int^{\tau}_0\int_{\Omega}\cA(t,x, \nabla u)\cdot \nabla T_{k}(u)\,dxdt.
\end{split}
\end{equation}

For the term on the right-hand side of equation~\eqref{eqnew}, we have
\begin{equation}\label{eqnew3}
\lim\limits_{l\to +\infty}\lim\limits_{\delta\to 0}\lim\limits_{n\to +\infty}\int^{\tau}_0\int_{\Omega}f_{n}\psi_{l}(u_n)(T_{k}(u))_{\delta}\,dxdt=\int^{\tau}_0\int_{\Omega}fT_{k}(u)\,dxdt.
\end{equation}

Thus, it follows from~\eqref{eqnew}--\eqref{eqnew3} that
\begin{equation}\label{eq:322}
\int^{\tau}_0\int_{\Omega}\partial_tu T_{k}(u)\,dxdt+\int^{\tau}_0\int_{\Omega}\cA(\tau,x, \nabla u)\cdot \nabla T_{k}(u)\,dxdt=\int^{\tau}_0\int_{\Omega}fT_{k}(u)\,dxdt.
\end{equation}

In addition, testing the approximate problem \eqref{eq:appro} by $T_{k}(u_{n})\chi_{(0,\tau)}$, we have
\begin{equation}\label{eq:321}
    \begin{split}
        & \lim_{n\to+\infty}\bigg(\int_{\Omega}\Theta_{k}(u_{n})(\tau)\,dx -\int_{\Omega}\Theta_{k}(u_{0n})\,dx+
\int^{\tau}_{0}\int_{\Omega}\cA(t,x,\nabla T_{k}(u_{n}))\cdot \nabla T_{k}(u_{n})\,dxdt\bigg)\\
&\quad =  \lim_{n\to+\infty}\int^{\tau}_{0}\int_{\Omega}f_{n}T_{k}(u_{n})\,dxdt=\int^{\tau}_{0}\int_{\Omega}fT_{k}(u)\,dxdt.
    \end{split}
\end{equation}

Combining with \eqref{eq:322} and~\eqref{eq:321}, we obtain
\begin{equation}\label{eqnewgoal}
  \lim_{n\to+\infty} \int^{\tau}_{0}\int_{\Omega}\cA(t,x,\nabla T_{k}(u_n))\cdot\nabla T_k(u_n)\,dxdt= \int^{\tau}_{0}\int_{\Omega}\cA(t,x,\nabla
 T_{k}(u))\cdot\nabla T_{k}(u)\,dxdt.
\end{equation}

Thus, Lemma \ref{zui} gives that
\begin{eqnarray*}
        \cA(t,x,\nabla T_{k}(u_n)) \cdot \nabla T_k(u_n)  \rightarrow \cA(t,x,\nabla T_{k}(u))\cdot \nabla T_{k}(u) ~~~~~\quad ~\quad ~~\mbox{strongly in }~~~~~~~L^{1}(\Omega_T).
    \end{eqnarray*}

\noindent\textbf{Step 5.} We shall prove that $u$ is a renormalized solution.

\textbf{Condition(R1)}. Firstly, we shall show that the  $u$ satisfies condition ($R1$). By the properties of truncations, for every $l\in \mathbb{N}$, we have
\begin{equation*}
 \begin{split}
 \nabla u_n=0 \ \ \ \text{a.e. in} \ \ \ \{(t,x) \in \Omega_T : |u_n|\in \{l,l+1\}\}.
\end{split}
\end{equation*}
From this and \eqref{eqstep3}, it follows that
  \begin{eqnarray}\label{R13}
\lim_{l  \to+ \infty}\lim\limits_{n\rightarrow +\infty} \int_{\{ l-1 \leq |u_n| \leq l+2\}} \cA(t,x,\nabla u_n)\cdot \nabla u_n \,dxdt=0.
 \end{eqnarray}
Now, define the function $G_l:\mathbb{R}\to \mathbb{R}$ by:
$$G_l(r)=\left\{
 \begin{array}{ll}
 1 & \ \ \ \mbox{if }  l\leq |r|\leq l+1, \\
 0 & \ \ \ \mbox{if }  |r|<l-1 \ \text{or} \ |r|>l+2,\\
 \text{affine} &  \ \ \ \text{otherwise}.
 \end{array}
\right.$$
 Using this, we can write
 \begin{equation*}
 \begin{split}
 \int_{\{l<|u|<l+1\}}\cA(t,x,\nabla u)\cdot \nabla u \,dxdt \leq \int^{T}_{0}\int_{\Omega}G_l(u)\cA(t,x,\nabla T_{l+2}(u))\cdot \nabla T_{l+2}(u)\,dxdt.
 \end{split}
\end{equation*}
Taking the limit in the inequality above gives
 \begin{equation}\label{R12}
 \begin{split}
0\leq &\lim\limits_{l\to +\infty} \int_{\{l<|u|<l+1\}}\cA(t,x,\nabla u)\cdot \nabla u \,dxdt \\
\leq &\lim\limits_{l\to +\infty} \int^{T}_{0}\int_{\Omega}G_l(u)\cA(t,x,\nabla T_{l+2}(u))\cdot \nabla T_{l+2}(u)\,dxdt.
 \end{split}
\end{equation}
Using \eqref{goal}, \eqref{R13} and the fact that $G_l$  is bounded and continuous, we get
 \begin{equation}\label{R11}
 \begin{split}
&\lim\limits_{l\to +\infty}\int^{T}_{0}\int_{\Omega}G_l(u)\cA(t,x,\nabla T_{l+2}(u))\cdot \nabla T_{l+2}(u)\,dxdt\\
 &\quad=\lim\limits_{l\to +\infty}\lim\limits_{n\to +\infty}\int^{T}_{0}\int_{\Omega}G_l(u)\cA(t,x,\nabla T_{l+2}(u_n))\cdot \nabla T_{l+2}(u_n)\,dxdt\\
 &\quad \leq \lim\limits_{l\to +\infty}\lim_{n\to +\infty}\int_{\{l-1<|u_n|<l+2\}}\cA(t,x,\nabla u_n)\cdot \nabla u_n \,dxdt \\
 &\quad=0.
 \end{split}
\end{equation}
From~\eqref{R12} and \eqref{R11}, we conclude that
 \begin{equation*}
 \lim\limits_{l\to +\infty}\int_{\{l<|u|<l+1\}}\cA(t,x,\nabla u)\cdot \nabla u \,dxdt=0,
\end{equation*}
 which gives condition (R1).

 \vskip .2cm

\textbf{Condition(R2)}.
Let $S\in W^{2,\infty}(\mathbb{R})$ such that $\text{supp} S'\in [-M,M]$ for some $M>0$. For every $ \phi \in C^{1}(\overline{\Omega}_{T})$ with $\phi(\cdot, T)=0$, taking $S'(u_n)\phi$ as a test function in \eqref{eq:appro}, we obtain
\begin{equation*}\label{eq:r}
    \begin{split}
        &\int^{T}_{0}\int_{\Omega} \frac{\partial S(u_n)}{\partial t}\phi\,dxdt+\int^{T}_{0}\int_{\Omega}\cA(t,x,\nabla u_n)\cdot\nabla\phi S'(u_n)+\cA(t,x,\nabla u_n)\cdot\nabla u_n S''(u_n)\phi\,dxdt\\
        &\quad=\int^{T}_{0}\int_{\Omega}f_n S'(u_n)\phi\,dxdt.
    \end{split}
\end{equation*}
Recalling the fact that $u_n\to u$ a.e. in $\Omega_T$, we obtain that
\begin{equation*}
    \int^{T}_{0}\int_{\Omega} \frac{\partial S(u_n)}{\partial t}\phi\,dxdt \to \int^{T}_{0}\int_{\Omega} \frac{\partial S(u)}{\partial t}\phi\,dxdt
\end{equation*}
and
\begin{equation*}
 \int^{T}_{0}\int_{\Omega}f_n S'(u_n)\phi\,dxdt\to \int^{T}_{0}\int_{\Omega}f S'(u_n)\phi\,dxdt
\end{equation*}
as $n\to +\infty$.

In addition, since $\supp S'\in [-M,M]$ for $M>0$, we know that
\[S'(u_n) \cA(t,x,\nabla u_n)= S'(u_n)\cA(t,x,\nabla T_{M}(u_n)).\]
Note that $S'(u_n)$ is bounded and satisfies 
\[S'(u_n)\to S'(u)~\quad ~\text{a.e. in}~\Omega_T.\]
Thus, it follows from
\[\cA(t,x,\nabla T_{M}(u_n))\rightharpoonup\cA(t,x,\nabla T_{M}(u))~\quad ~\text{weakly~in~} L^1(\Omega_T)\]
that
\[S'(u_n)\cA(t,x,\nabla T_{M}(u_n))\rightharpoonup S'(u)\cA(t,x,\nabla T_{M}(u))~\text{weakly~in~} L^1(\Omega_T),\]
which implies that
\begin{equation*}
    \int^{T}_{0}\int_{\Omega}\cA(t,x,\nabla u_n) S'(u_n)\cdot\nabla\phi\,dxdt \to \int^{T}_{0}\int_{\Omega}\cA(t,x,\nabla u) S'(u)\cdot \nabla\phi\,dxdt
\end{equation*}
as $n\to +\infty$. Moreover,  it follows from \eqref{goal}  that
\begin{equation*}
   \int^{T}_{0}\int_{\Omega}\cA(t,x,\nabla u_n)\cdot\nabla u_n S''(u_n)\phi\,dxdt\to   \int^{T}_{0}\int_{\Omega}\cA(t,x,\nabla u)\cdot\nabla u S''(u)\phi\,dxdt
\end{equation*}
as $n\to +\infty$.  Thus,  we conclude that
\begin{equation*}\label{eq:r2}
    \begin{split}
        &\int^{T}_{0}\int_{\Omega} \frac{\partial S(u)}{\partial t}\phi\,dxdt+\int^{T}_{0}\int_{\Omega}\left(\cA(t,x,\nabla u)\cdot\nabla\phi S'(u)+\cA(t,x,\nabla u)\cdot\nabla u S''(u)\phi\right)\,dxdt\\
        &\quad=\int^{T}_{0}\int_{\Omega}f S'(u)\phi\,dxdt.
    \end{split}
\end{equation*}

We have proven the existence of the renormalized solution. Now, we will proceed to demonstrate its uniqueness.

\vskip .2cm

\textbf{(2). Uniqueness of renormalized solutions.}

Let $u$ and $v$ be two renormalized solutions to problem \eqref{eq:main}. For $\sigma >0$, define the function $S_{\sigma}$ as follows:
\begin{align*}
			S_\sigma(r):=
			\begin{cases}
				r&\text{if }|r|<\sigma,\\
				\left( \sigma+\frac{1}{2}\right)-\frac{1}{2}\left( r-(\sigma+1)\right)^2  &\text{if }\sigma\le r\le\sigma+1,\\
				-\left( \sigma+\frac{1}{2}\right)+\frac{1}{2}\left( r+(\sigma+1)\right)^2  &\text{if }-\sigma-1\le r\le-\sigma,\\
				\left( \sigma+\frac{1}{2}\right)  &\text{if }r>\sigma+1,\\
				-\left( \sigma+\frac{1}{2}\right)  &\text{if }r<-\sigma-1.
			\end{cases}
		\end{align*}
  It is obvious that
  \begin{align*}
			S_\sigma'(r)=
			\begin{cases}
				1  &\text{if }|r|< \sigma,\\
				\sigma+1-|r| &\text{if }\sigma\le |r|\le \sigma+1,\\
				0 &\text{if }|r|>\sigma+1,
			\end{cases}
		\end{align*}
$S_\sigma\in W^{2,\infty}(\mathbb{R})$ with $\text{supp}\,S_\sigma'\subset [-\sigma-1,\sigma+1]$ and $\text{supp}\,S_\sigma''\subset [\sigma,\sigma+1]\cup[-\sigma-1,-\sigma]$.

We take $S=S_{\sigma}$ to obtain
\begin{equation}
    \int^{T}_{0}\int_{\Omega} \frac{\partial S_{\sigma}(u)}{\partial t}\phi \,dxdt +\int^{T}_{0}\int_{\Omega}\cA(t,x,\nabla u)\cdot
    \nabla (S'(u)\phi)\,dxdt=\int^{T}_{0}\int_{\Omega}fS'_{\sigma}(u)\phi\,dxdt
\end{equation}
and
\begin{equation}
    \int^{T}_{0}\int_{\Omega} \frac{\partial S_{\sigma}(v)}{\partial t}\phi\,dxdt +\int^{T}_{0}\int_{\Omega}\cA(t,x,\nabla v)\cdot
    \nabla (S'(v)\phi)\,dxdt=\int^{T}_{0}\int_{\Omega}fS'_{\sigma}(v)\phi\,dxdt.
\end{equation}
For fixed $k>0$, taking $\phi=T_{k}(S_{\sigma}(u)-S_{\sigma}(v))$ and subtracting the two above equations, we deduce that
\begin{equation*}
    \begin{split}
        &\int^{T}_{0}\int_{\Omega} \left(\frac{\partial S_{\sigma}(u)}{\partial t}-\frac{\partial S_{\sigma}(v)}{\partial t}\right)T_{k}(S_{\sigma}(u)-S_{\sigma}(v))\,dxdt \\
        &\quad +\int^{T}_{0}\int_{\Omega}\cA(t,x,\nabla u)\cdot\nabla u
        S''(u)T _{k}(S_{\sigma}(u) -S_{\sigma}(v))\\
        &\quad \quad \quad \quad \quad-\cA(t,x,\nabla v)\cdot\nabla v S''(v)T_{k}(S_{\sigma}(u)-S_{\sigma}(v))\,dxdt\\
        &\quad+\int^{T}_{0}\int_{\Omega}S'(u)\cA(t,x,\nabla u)\cdot\nabla T_{k}(S_{\sigma}(u)-S_{\sigma}(v)) \\
        &\quad \quad \quad \quad \quad -S'(v)\cA(t,x,\nabla v)\cdot\nabla T_{k}(S_{\sigma}(u)-S_{\sigma}(v))\,dxdt\\
        &=\int^{T}_{0}\int_{\Omega}(fS'_{\sigma}(u)-fS'_{\sigma}(v))T_{k}(S_{\sigma}(u)-S_{\sigma}(v))\,dxdt.
    \end{split}
\end{equation*}
We set
\begin{align*}
    &J_0:=\int^{T}_{0}\int_{\Omega}\left(\frac{\partial S_{\sigma}(u)}{\partial t}-\frac{\partial S_{\sigma}(v)}{\partial t}\right)T_{k}(S_{\sigma}(u)-S_{\sigma}(v))\,dxdt, \\
    & J_1:=\int^{T}_{0}\int_{\Omega}\cA(t,x,\nabla u)\cdot\nabla u
        S_{\sigma}''(u)T _{k}(S_{\sigma}(u) -S_{\sigma}(v))\\
        &\quad \quad \quad \quad \quad-\cA(t,x,\nabla v)\cdot\nabla v S_{\sigma}''(v)T_{k}(S_{\sigma}(u)-S_{\sigma}(v))\,dxdt,\\
   &J_2:=\int^{T}_{0}\int_{\Omega}S_{\sigma}'(u)\cA(t,x,\nabla u)\cdot\nabla T_{k}(S_{\sigma}(u) -S_{\sigma}(v))\\
   &\quad \quad \quad \quad \quad-S_{\sigma}'(v)\cA(t,x,\nabla v)\cdot\nabla T_{k}(S_{\sigma}(u)-S_{\sigma}(v))\,dxdt,\\
   &J_3:=\int^{T}_{0}\int_{\Omega}(fS'_{\sigma}(u)-fS'_{\sigma}(v))T_{k}(S_{\sigma}(u)-S_{\sigma}(v))\,dxdt.
\end{align*}
Thus, we have
\begin{equation*}
    J_0+J_1+J_2=J_3.
\end{equation*}
Next, we shall estimate $J_i(i=0,1,2,3)$ one by one.

 \textit{Estimates of $J_0$.} Due to the same initial condition for $u$ and $v$, and the properties of $\Theta_k$, we get
 \begin{equation*}
     J_0=\int_\Omega\Theta_k\left( S_\sigma(u)-S_\sigma(v)\right)(T)dx\ge 0.
 \end{equation*}

\textit{Estimates of $J_1$.} It follows from the definition of $S_{\sigma}$ that
\begin{equation*}
|J_1|\leq Ck\left(\int_{\{\sigma<|u|<\sigma+1\}} \cA(t,x,\nabla u)\cdot\nabla u\,dxdt+\int_{\{\sigma<|v|<\sigma+1\}} \cA(t,x,\nabla v)\cdot\nabla v\,dxdt\right).
\end{equation*}
According to the definition of renormalized solutions, we have
\begin{equation*}
    \lim_{\sigma\to+\infty} J_1=0.
\end{equation*}

\textit{Estimates of $J_2$.}  We write
\begin{equation*}
    \begin{split}
        J_2&:=\int^{T}_{0}\int_{\Omega} \cA(t,x,\nabla u) (S'_{\sigma}(u)-1)\cdot\nabla T_{k}(S_{\sigma}(u)-S_{\sigma}(v))\,dxdt\\
        &\quad +\int^{T}_{0}\int_{\Omega}(\cA(t,x,\nabla u)-\cA(t,x,\nabla v))\cdot\nabla T_{k}(S_{\sigma}(u)-S_{\sigma}(v))\,dxdt\\
        &\quad+\int^{T}_{0}\int_{\Omega} \cA(t,x,\nabla v) (1-S'_{\sigma}(v))\cdot\nabla T_{k}(S_{\sigma}(u)-S_{\sigma}(v))\,dxdt\\
        &:= J^1_2+J^2_2+J^2_3.
    \end{split}
\end{equation*}
Setting $\sigma >k$, we obtain
\begin{equation*}
    J^2_2 \geq \int_{\{|u-v|\leq k\cap |u|\leq k,|v|\leq k\}}(\cA(t,x,\nabla u)-\cA(t,x,\nabla v)\cdot\nabla (u-v)\,dxdt.
\end{equation*}
Moreover, since meas($\{u=\infty\}=0$) and meas($\{v=\infty\}=0$), we deduce that
\begin{equation*}
\begin{split}
     &\lim_{\sigma\to+\infty}(|J^1_2|+|J^3_2|)\\
     &\quad \leq \lim_{\sigma\to+\infty} \int_{\{\sigma\leq |u|\}} \cA(t,x,\nabla u) \cdot\nabla T_{k}(S_{\sigma}(u)-S_{\sigma}(v))\,dxdt\\
      &\qquad+ \lim_{\sigma\to+\infty} \int_{\{\sigma\leq |v|\}} \cA(t,x,\nabla v)\cdot \nabla T_{k}(S_{\sigma}(u)-S_{\sigma}(v))\,dxdt\\
      &\quad =0.
\end{split}
\end{equation*}
Thus, we conclude that
\begin{equation*}
    \lim_{\sigma\to+\infty} J_2\geq \int_{\{|u-v|\leq k\cap |u|\leq k,|v|\leq k\}}(\cA(t,x,\nabla u)-\cA(t,x,\nabla v)\cdot \nabla (u-v)\,dxdt.
\end{equation*}

\textit{Estimates of $J_3$.} It is easy to check that
\begin{equation*}
    |J_3|\leq k\int_{\{\text{max}\{|u|,|v|\}>\sigma\}}|f|\,dxdt,
\end{equation*}
which implies that
\begin{equation*}
    \lim_{\sigma\to+\infty} |J_3|=0.
\end{equation*}

To sum up, by sending $\sigma\to+\infty$, we obtain
\begin{equation*}
    \begin{split}
        0&\leq \int_{\{|u|\leq \frac{k}{2}, |v|\leq\frac{k}{2}\}}(\cA(t,x,\nabla u)-\cA(t,x,\nabla v)\cdot\nabla (u-v)\,dxdt\\
        &\leq 0,
    \end{split}
\end{equation*}
which implies that $\nabla u=\nabla v$ a.e. on the set $\{|u|\leq \frac{k}{2}, |v|\leq\frac{k}{2}\}$. Since $k$ is arbitrary, we have $\nabla u=\nabla v$ a.e. in $\Omega_T$. Thus, from the Poincar\'{e} inequality, we have $u=v$ a.e. in $\Omega_T$. \hfill $\Box.$

\vskip .2cm

We have proven the existence and uniqueness of the renormalized solution to equation~\eqref{eq:main}. Next, we will demonstrate that this renormalized solution is also its entropy solution and prove the uniqueness of the entropy solution. Consequently, we obtain the equivalence between the entropy and renormalized solutions for this equation.

\vskip .2cm

\noindent\textbf{Proof of Theorem \ref{thm:entropy}}. First, we will prove that the renormalized solution of equation \eqref{eq:main} is also its entropy solution.

\textbf{(1). Existence of entropy solutions.}

 Now we choose $T_{k}(u_{n}-\phi)$ as a test function in \eqref{eq:appro} for $k\in \mathbb{N}^{+}$ and  $\phi\in C^{1}(\overline{\Omega}_{T})$ with $\phi|_{\Sigma}=0$. Set $L=k+\|\phi\|_{L^{\infty}(\Omega_{T})}$, we deduce
\begin{eqnarray*}
\int^{T}_{0}\int_{\Omega}\cA(t,x,\nabla u_{n})\cdot \nabla T_{k}(u_{n}-\phi)\,dxdt =\int^{T}_{0}\int_{\Omega}\cA(t,x,\nabla T_{L}(u_{n}))\cdot \nabla T_{k}(T_{L}(u_{n})-\phi)\,dxdt
\end{eqnarray*}
and
\begin{eqnarray*}
&&\int^{T}_{0}\langle \partial_t u_{n},T_{k}(u_{n}-\phi)\rangle \,dt+\int^{T}_{0}\int_{\Omega}\cA\left(t,x,\nabla T_{L}(u_{n})\right)\cdot \nabla T_{k}\left(T_{L}(u_{n})-\phi\right)\,dxdt \nonumber\\
&& \quad =\int^{T}_{0}\int_{\Omega}f_{n}T_{k}(u_{n}-\phi)\,dxdt.
\end{eqnarray*}
Since $\partial_t u_{n}=\partial_t(u_{n}-\phi)+\partial_t\phi$, we have
\begin{eqnarray*}
\int^{T}_{0}\langle \partial_t u_{n},T_{k}(u_{n}-\phi)\rangle \,dt &=& \int_{\Omega}\Theta_{k}(u_{n}-\phi)(T)dx- \int_{\Omega}\Theta_{k}(u_{n}-\phi)(0)\,dx  \\
& &\quad +\int^{T}_{0}\langle \partial_t\phi,T_{k}(u_{n}-\phi)\rangle dt,
\end{eqnarray*}
which yields that
\begin{equation}\label{316}
\begin{split}
&\int_{\Omega}\Theta_{k}(u_{n}-\phi)(T)\,dx -\int_{\Omega}\Theta_{k}(u_{n}-\phi)(0)\,dx +\int^{T}_{0}\langle\phi_{t},T_{k}(u_{n}-\phi)\rangle\,dt\\
& \quad +\int^{T}_{0}\int_{\Omega}\cA(t,x,\nabla T_{L}(u_{n}))\cdot \nabla T_{k}(T_{L}(u_{n})-\phi)\,dxdt =\int^{T}_{0}\int_{\Omega}f_{n}T_{k}(u_{n}-\phi)\,dxdt.
\end{split}\end{equation}
Recalling $u_{n}$ converges to $u$ in $C([0,T];L^{1}(\Omega))$, we have $u_{n}(t)\rightarrow u(t)$ in $L^{1}(\Omega)$  as $n  \to+ \infty$ for all $t \le T$. Since $\Theta_{k}$ is Lipschitz continuous, we get
\begin{eqnarray*}
\int_{\Omega}\Theta_{k}(u_{n}-\phi)(T)dx \rightarrow \int_{\Omega}\Theta_{k}(u-\phi)(T)dx
\end{eqnarray*}
and
\begin{eqnarray*}
\int_{\Omega}\Theta_{k}(u_{n}-\phi)(0)dx \rightarrow \int_{\Omega}\Theta_{k}(u_{0}-\phi(0))dx
\end{eqnarray*}
as $n\rightarrow +\infty$.

The fourth term on the left-hand side of (\ref{316}) can be written as
\begin{equation}\label{315}
\begin{split}
&\int^{T}_{0}\int_{\Omega}\cA(t,x,\nabla T_{L}(u_{n})) \cdot \nabla T_{k}(T_{L}(u_{n})-\phi)\,dxdt\\
&\quad=\int_{\{|T_{L}(u_{n})-\phi|\le k\}}\cA(t,x,\nabla T_{L}(u_{n})) \cdot \nabla T_{L}(u_{n})\,dxdt\\
&\quad \quad -\int_{\{|T_{L}(u_{n})-\phi|\le k\}}\cA(t,x,\nabla T_{L}(u_{n})) \cdot  \nabla \phi \,dxdt\\
&\quad =:I_{1}+I_{2},
\end{split}
\end{equation}
where
\begin{eqnarray*}
I_{1}:=\int_{\{|T_{L}(u_{n})-\phi|\le k\}}\cA(t,x,\nabla T_{L}(u_{n})) \cdot \nabla T_{L}(u_{n})\,dxdt
\end{eqnarray*}
and
\begin{eqnarray*}
I_{2}:=-\int_{\{|T_{L}(u_{n})-\phi|\le k\}}\cA(t,x,\nabla T_{L}(u_{n})) \cdot  \nabla \phi\,dxdt.
\end{eqnarray*}
\noindent\emph{Estimate of $I_{1}$}.
 We derive from \eqref{goal} that
\begin{equation}\label{2222}
\begin{split}
&\int_{\{|T_{L}(u)-\phi|\le k\}}\cA(t,x,\nabla T_{L}(u))\cdot\nabla T_{L}(u)\,dxdt \\
&\quad =\displaystyle \lim_{n \to+ \infty}\int_{\{|T_{L}(u_{n})-\phi|\le k\}}\cA(t,x,\nabla T_{L}(u_{n})) \cdot \nabla T_{L}(u_{n})\,dxdt \\
&\quad =:\displaystyle\lim_{n \to+ \infty} I_{1}.
\end{split}
\end{equation}
\noindent\emph{Estimate of $I_{2}$}. For convenience, we define
$$ \eta_{n}:= -\cA(t,x,\nabla T_{L}(u_{n})),\quad E_{n}:=\{ (t,x)\in\Omega_{T}:|T_{L}(u_{n})-\phi|\leq k\} $$
and
$$E:=\{ (t,x)\in \Omega_{T}:|T_{L}(u)-\phi|\leq k\}.$$
We can write
\begin{align*}
I_{2}=\int_{E_{n}}\eta_{n}\cdot \nabla \phi \,dxdt= \int_{E}\eta_{n}\cdot \nabla \phi dz +\int_{E_{n}\setminus E}\eta_{n}\cdot \nabla \phi \,dxdt:= I_{21}+I_{22}.
\end{align*}

Recalling the fact that  for  $n \to+ \infty$
\begin{eqnarray}\label{2.9000}
\cA(t,x,\nabla T_{L}(u_{n})) \rightharpoonup \cA(t,x,\nabla T_{L}(u))~ ~\qquad~~ ~\mbox{weakly  in}~~~~~L^{1}(\Omega_{T}),
\end{eqnarray}
we have
\begin{align*}
\lim_{n  \to+ \infty}I_{21}=-\int_{\{|T_{L}(u)-\phi|\leq k\}}\cA(t,x,\nabla T_{L}(u))\cdot \nabla \phi \,dxdt.
\end{align*}
Moreover, since $M^{*}(t,x,\eta)$ satisfies the condition $\displaystyle \lim_{|\eta|\rightarrow +\infty} \frac{M^{*}(t,x, |\eta|)}{|\eta|}=+\infty$. Then for every $\varepsilon >0$, there exists a constant $\Lambda>0$ such that
\begin{align*}
|\eta| \leq \varepsilon M^{*}(t,x,|\eta|)~~~~~~\quad ~~~\mbox{for~~~~all}  ~~~~~ |\eta|>\Lambda.
\end{align*}
It follows from (\ref{3.5}) that
\begin{align*}
|I_{22}| &\leq \|\nabla\phi\|_{L^{\infty}(\Omega_{T})}\int^{T}_{0}\int_{\Omega}|\eta_{n}| \chi_{E_{n}\setminus E} \,dxdt \\
&= C\left(\int_{\{|\eta_{n}|\leq \Lambda\}}|\eta_{n}|\chi_{E_{n}\setminus E}\,dxdt+\int_{\{|\eta_{n}|> \Lambda \}}|\eta_{n}|\chi_{E_{n}\setminus E}\,dxdt\right)\\
&\leq C\left(\Lambda\mbox{meas}(E_{n}\setminus E)+\varepsilon \int^T_0\int_{\Omega}M^{*}(t,x, |\eta_{n}|)\,dxdt\right)\\
&\leq C\Lambda \mbox{meas}(E_{n}\setminus E)+C\varepsilon.
\end{align*}
By the arbitrariness of $\varepsilon$, we get
$$ \lim_{n\rightarrow+\infty}|I_{22}|=0. $$
Therefore, we have
\begin{eqnarray}\label{21000}
 \lim_{n\rightarrow+\infty} I_{2} =-\int_{\{|T_{L}(u)-\phi|\le k\}}\cA(t,x,\nabla T_{L}(u))\cdot\nabla \phi\,dxdt.
 \end{eqnarray}
According to (\ref{2222}) and (\ref{21000}), we obtain
\begin{eqnarray*}
\int^{T}_{0}\int_{\Omega}\cA(t,x,\nabla u)\cdot \nabla T_{k}(u-\phi)\,dxdt &=& \int^{T}_{0}\int_{\Omega}\cA(t,x,\nabla T_{L}(u))\cdot \nabla T_{k}(T_{L}(u)-\phi)\,dxdt\nonumber\\
&=& \int_{\{|T_{L}(u)-\phi|<k\}}\cA(t,x,\nabla T_{L}(u))\cdot\nabla T_{L}(u)\,dxdt \nonumber\\
& &-\int_{\{|T_{L}(u)-\phi|<k\}}\cA(t,x,\nabla T_{L}(u))\cdot\nabla \phi\,dxdt\nonumber \\
&=& \lim_{n\rightarrow+\infty}(I_{1}+I_{2}).
\end{eqnarray*}
Using the strong convergence of $f_{n}$, (\ref{3}) and the Lebesgue dominated convergence theorem, we can pass to the limits as $n\rightarrow +\infty$ in the other term of (\ref{316}) to conclude
\begin{eqnarray}
&&\int_{\Omega}\Theta_{k}(u-\phi)(T)dx -\int_{\Omega}\Theta_{k}(u_{0}-\phi(0))dx +\int^{T}_{0}\langle \partial_t\phi,T_{k}(u-\phi)\rangle \,dt \nonumber \\
& & \quad +\int^{T}_{0}\int_{\Omega}\cA(t,x,\nabla u)\cdot\nabla T_{k}(u-\phi)\,dxdt = \int^{T}_{0}\int_{\Omega}fT_{k}(u-\phi)\,dxdt
\end{eqnarray}
for all $k>0$ and $\phi\in C^{1}(\overline{\Omega}_{T})$ with $\phi|_{\Sigma}=0$. Hence, our solution $u$ satisfies condition $(E2)$. Since $u_{n}$ is the distributional solution of problem \eqref{eq:appro}, then its limit $u$ satisfies condition $(E1)$ naturally. This completes the proof.

\vskip.2cm

We have established the existence of the entropy solution. Now, we further show its uniqueness.

\vskip.2cm

\textbf{(2). Uniqueness of entropy solutions.}

Suppose that $v$ is another entropy solution of problem \eqref{eq:main}, we will show that $u=v$ a.e. in $\Omega_{T}$. For $\sigma>0, 0<\varepsilon\le  1$, define the function $S_{\sigma,\varepsilon}$ in $W^{2,\infty}(\mathbb{R})$ by
\begin{align*}
  S_{\sigma,\varepsilon}(r):= \left\{
  \begin{array}{ccc}
r \quad \quad \quad  &\mbox{if}& |r|\le \sigma,\\
\left(\sigma+\frac{\varepsilon}{2}\right)- \frac{r}{|r|}\frac{1}{2\varepsilon}\left(r- \frac{r}{|r|}(\sigma +\varepsilon)\right)^{2} \quad \quad &\mbox{if}& \sigma< |r|<\sigma+\varepsilon, \\
\frac{r}{|r|}\left(\sigma+\frac{\varepsilon}{2}\right) \quad \quad \quad &\mbox{if}& |r|\ge \sigma+\varepsilon.
 \end{array}
   \right.
\end{align*}
Clearly,
\begin{align*}
S'_{\sigma,\varepsilon}(r)=\left\{
\begin{array}{ccc}
1 \quad \quad \quad &\mbox{if}& |r|\le \sigma,\\
\frac{1}{\varepsilon}\left( \sigma+\varepsilon -|r|\right) \quad \quad \quad &\mbox{if}& \sigma < |r| <\sigma +\varepsilon, \\
0 \quad \quad \quad &\mbox{if}& |r|\ge \sigma+\varepsilon.
\end{array}
\right.
\end{align*}

\noindent Choosing $\phi=S_{\sigma,\varepsilon}(u_{n})\chi_{(0,\tau)}$ as a test function in (\ref{f5}) for entropy solution $v$, we have
\begin{eqnarray}\label{319}
H_{0}+H_{1}+H_{2} = H_{3},
\end{eqnarray}
where
\begin{eqnarray*}
&&H_{0}:=\int_{\Omega}\Theta_{k}(v- S_{\sigma,\varepsilon}(u_{n}))(\tau)dx -\int_{\Omega}\Theta_{k}(u_{0}-S_{\sigma,\varepsilon}(u_{0n}))\,dx,\nonumber\\
&&H_{1}:=\int^{\tau}_{0}\langle \partial_t u_{n}, S'_{\sigma,\varepsilon}(u_{n})T_{k}\left(v-S_{\sigma,\varepsilon}(u_{n})\right)\rangle \,dt \nonumber, \\
&&H_{2}:=\int^{\tau}_{0}\int_{\Omega}\cA(t,x,\nabla v)\cdot \nabla T_{k}\left(v-S_{\sigma,\varepsilon}(u_{n})\right)\,dxdt, \nonumber \\
&&H_{3}:=\int^{\tau}_{0}\int_{\Omega}fT_{k}\left(v-S_{\sigma,\varepsilon}(u_{n})\right)\,dxdt.
\end{eqnarray*}
Taking $S'_{\sigma,\varepsilon}(u_{n})T_{k}\left(v-S_{\sigma,\varepsilon}(u_{n})\right)\chi_{(0,\tau)}$ as a test function for problem \eqref{eq:appro} we get
\begin{eqnarray}\label{320}
&&H_{1}= \int^{\tau}_{0}\langle \partial_t u_{n},S'_{\sigma,\varepsilon}(u_{n})T_{k}(v-S_{\sigma,\varepsilon}(u_{n}))\rangle\,dt \nonumber \\
&&\quad =\int^{\tau}_{0}\int_{\Omega}f_{n}S'_{\sigma,\varepsilon}(u_{n}) T_{k}\left( v-S_{\sigma,\varepsilon}(u_{n})\right)\,dxdt\nonumber\\
&& \quad \quad -\int^{\tau}_{0}\int_{\Omega}S''_{\sigma,\varepsilon}(u_{n})T_{k}\left(v-S_{\sigma,\varepsilon}(u_{n})\right)\cA(t,x,\nabla u_{n})\cdot\nabla u_{n}\,dxdt \nonumber \\
&& \quad \quad -\int^{\tau}_{0}\int_{\Omega}S'_{\sigma,\varepsilon}(u_{n})\cA(t,x,\nabla u_{n})\cdot \nabla T_{k}\left(v-S_{\sigma,\varepsilon}(u_{n})\right)\,dxdt \nonumber\\
&&\quad =: H_{11}+H_{12}+H_{13}.
\end{eqnarray}
Combining with (\ref{319}) and (\ref{320}), we have
\begin{eqnarray*}
&&\int_{\Omega}\Theta_{k}\left(v-S_{\sigma,\varepsilon}(u_{n})\right)(\tau)\,dx -\int_{\Omega}\Theta_{k}\left(u_{0}-S_{\sigma,\varepsilon}(u_{0n})\right)\,dx \nonumber \\
&& \quad-\int^{\tau}_{0}\int_{\Omega}S'_{\sigma,\varepsilon}(u_{n})\cA(t,x,\nabla u_{n})\cdot \nabla T_{k}\left(v-S_{\sigma,\varepsilon}(u_{n})\right)\,dxdt\nonumber \\
&& \quad +\int^{\tau}_{0}\int_{\Omega}\cA(t,x,\nabla v)\cdot \nabla T_{k}\left(v-S_{\sigma,\varepsilon}(u_{n})\right)\,dxdt \nonumber \\
&&\le \int^{\tau}_{0}\int_{\Omega}fT_{k}\left(v-S_{\sigma,\varepsilon}(u_{n})\right)\,dxdt \nonumber \\
&&\quad -\int^{\tau}_{0}\int_{\Omega}f_{n}S'_{\sigma,\varepsilon}(u_{n}) T_{k}\left( v-S_{\sigma,\varepsilon}(u_{n})\right)\,dxdt\\
&&\quad+\left| \int^{\tau}_{0}\int_{\Omega}S''_{\sigma,\varepsilon}(u_{n})T_{k}\left(v-S_{\sigma,\varepsilon}(u_{n})\right)\cA(t,x,\nabla u_{n})\cdot\nabla u_{n}\,dxdt \nonumber\right|.
\end{eqnarray*}
That is,
\begin{eqnarray}\label{791}
H_{0}+H_{13}+H_{2}\le H_{3}-H_{11}+|H_{12}|.
\end{eqnarray}
 We will pass to the limits as $\varepsilon \rightarrow 0$, $n  \to+ \infty$ and $\sigma  \to+ \infty$ successively.

 We begin with $\varepsilon \rightarrow 0$. Since $\left|\Theta_{k}\left(v-S_{\sigma,\varepsilon}(u_{n})\right)(\tau)\right|\le k\left(|v|+|T_{\sigma+1}(u_{n})|\right)(\tau)$, $S'_{\sigma,\varepsilon}(r)\le T'_{\sigma+1}(r)$ and $\left|\nabla T_{k}\left(v-S_{\sigma,\varepsilon}(u_{n})\right)\right| \le |\nabla T_{\sigma+k+1}(v)|+|\nabla T_{\sigma+1}(u_{n})|$, the three terms of the left-hand side and the two terms of the right-hand side in (\ref{791}) pass to the limit for $\varepsilon \rightarrow 0$ by the Lebesgue dominated convergence theorem. Now we estimate $|H_{12}|$. Let $R_{\sigma,\varepsilon}$ be an even function such that $R_{\sigma,\varepsilon}(r)=r-S_{\sigma,\varepsilon}(r)$ for $r \ge 0$. Then we choose $R'_{\sigma,\varepsilon}(u_{n})\chi_{(0,\tau)}$ as a test function in \eqref{eq:appro} to have
 \begin{eqnarray*}
 && \int_{\Omega}R_{\sigma,\varepsilon}(u_{n})(\tau)\,dx- \int_{\Omega}R_{\sigma,\varepsilon}(u_{0n})\,dx\\
 && \quad + \int^{\tau}_{0}\int_{\Omega}R''_{\sigma,\varepsilon}(u_{n})\cA(t, x,\nabla u_{n})\cdot \nabla u_{n}\,dxdt =\int^\tau_{0}\int_{\Omega}f_{n}R'_{\sigma,\varepsilon}(u_{n})\,dxdt.
 \end{eqnarray*}

\noindent Since $R_{\sigma,\varepsilon}(r)\ge 0$, $R_{\sigma,\varepsilon}(r)\le |r|$ on the set $\{|r|>\sigma\}$ and $|S''_{\sigma,\varepsilon}(r)|=R''_{\sigma,\varepsilon}(r)$, we obtain that
 \begin{eqnarray*}
 \int^{\tau}_{0}\int_{\Omega}|S''_{\sigma,\varepsilon}(u_{n})|\cA(t,x,\nabla u_{n})\cdot \nabla u_{n}\,dxdt &=& \int^{\tau}_{0}\int_{\Omega}R''_{\sigma,\varepsilon}(u_{n})\cA(t,x,\nabla u_{n})\cdot \nabla u_{n}\,dxdt \\
 &\le& \int_{\{|u_{n}|>\sigma\}}|f_{n}|\,dxdt+\int_{\{|u_{0n}|>\sigma\}}|u_{0n}|\,dx.
 \end{eqnarray*}
Therefore,
\begin{eqnarray*}
|H_{12}|\le k\left( \int_{\{|u_{n}|>\sigma\}}|f_{n}| \,dxdt+ \int_{\{|u_{0n}|>\sigma\}}|u_{0n}|dx\right).
\end{eqnarray*}
 Recalling that $\cA(t,x,0)=0$, $T'_{\sigma}(u_{n})\cA(t,x,\nabla u_{n})=\cA(t,x,\nabla T_{\sigma}(u_{n}))$,  and letting $\varepsilon \rightarrow 0$ in (\ref{791}), we obtain
\begin{eqnarray*}
&&\int_{\Omega} \Theta_{k}\left( v-T_{\sigma}(u_{n})\right)(\tau)dx -\int_{\Omega}\Theta_{k}\left(u_{0}-T_{\sigma}(u_{0n})\right)dx\\
&&\quad \quad +\int^{\tau}_{0}\int_{\Omega}\left(\cA(t,x,\nabla v)-\cA(t,x,\nabla T_{\sigma}(u_{n}))\right)\cdot \nabla T_{k}(v-T_{\sigma}(u_{n}))\,dxdt \\
&&\quad \le \int^{\tau}_{0}\int_{\Omega}\left(f-f_{n}T'_{\sigma}(u_{n})\right)T_{k}(v-T_{\sigma}(u_{n}))\,dxdt \\
&& \quad \quad +k \left(\int_{\{|u_{n}|>\sigma\}}|f_{n}|\,dxdt+\int_{\{|u_{0n}|>\sigma\}}|u_{0n}|dx\right).
\end{eqnarray*}

Next, we take $n\rightarrow +\infty$. Due to the fact that $\nabla T_{\sigma}(u_{n})\rightarrow \nabla T_{\sigma}(u)$ a.e. in $\Omega_{T}$ as $n\rightarrow +\infty$, Fatou's lemma and the Lebesgue dominated convergence theorem, sending $n\rightarrow +\infty$ in the above inequality we have
\begin{eqnarray*}
&&\int_{\Omega} \Theta_{k}\left( v-T_{\sigma}(u)\right)(\tau)\,dx -\int_{\Omega}\Theta_{k}\left(u_{0}-T_{\sigma}(u_{0})\right)\,dx\\
&&\quad \quad +\int^{\tau}_{0}\int_{\Omega}\left(\cA(t,x,\nabla v)-\cA(t,x,\nabla T_{\sigma}(u))\right)\cdot \nabla T_{k}(v-T_{\sigma}(u))\,dxdt \\
&&\quad \le \int^{\tau}_{0}\int_{\Omega}f\left(1-T'_{\sigma}(u)\right)T_{k}(v-T_{\sigma}(u))\,dxdt \\
&& \quad \quad +k \left(\int_{\{|u|>\sigma\}}|f|\,dxdt +\int_{\{|u_{0}|>\sigma\}}|u_{0}|dx\right).
\end{eqnarray*}

Now we let $\sigma  \to+ \infty$. Since
\begin{eqnarray*}
\left|\Theta_{k}(v-T_{\sigma}(u))(\tau)\right|\le k\left(|v(\tau)|+|u(\tau)|\right) \quad\mbox{and} \quad \left|\Theta_{k}\left(u_{0}-T_{\sigma}(u_{0})\right)\right|\le k|u_{0}|,
\end{eqnarray*}
by the Lebesgue dominated convergence theorem, we obtain
\begin{eqnarray*}
\int_{\Omega}\Theta_{k}(v-T_{\sigma}(u))(\tau)\,dx\rightarrow \int_{\Omega}\Theta_{k}(v-u)(\tau)\,dx, \quad \quad \quad \int_{\Omega}\Theta_{k}\left(u_{0}-T_{\sigma}(u_{0})\right)dx \rightarrow 0
\end{eqnarray*}
and
\begin{eqnarray*}
\int_{\{|u|>\sigma\}}|f|\,dxdt +\int_{\{|u_{0}|>\sigma\}}|u_{0}|\,dx \rightarrow 0.
\end{eqnarray*}

 Therefore, we conclude that
\begin{eqnarray*}
\int_{\Omega}\Theta_{k}(v-u)(\tau)dx+\int^{\tau}_{0}\int_{\Omega}\left(\cA(t,x,\nabla v)-\cA(t,x,\nabla u)\right)\cdot \nabla T_{k}(v-u)\,dxdt \le 0,
\end{eqnarray*}
which implies that
\begin{eqnarray*}
\int_{\Omega}\Theta_{k}(v-u)(\tau)dx +\int_{\{|u|\le \frac{k}{2},|v|\le \frac{k}{2}\}}\left(\cA(t,x,\nabla v)-\cA(t,x,\nabla u)\right)\cdot \nabla (v-u)\,dxdt \le 0.
\end{eqnarray*}
Using the nonnegativity of the two terms in the above inequality, we conclude that $u=v$ a.e. in $\Omega_{T}$. Therefore, we obtain the uniqueness of entropy solutions. \hfill $\Box.$

\subsection*{Acknowledgment}
 The first author was supported by the Postdoctoral Fellowship Program of CPSF under Grant Number: GZC20242215. The second author was supported by the National Natural Science Foundation of China (Nos. 12071098, 12471128).

\subsection*{Conflict of interest} The authors declare that there is no conflict of interest. We also declare that this manuscript has no associated data.

\subsection*{Data availability} Data sharing is not applicable to this article as no datasets were generated or analysed during the current study.

\end{document}